\documentclass[11pt]{amsart}
\usepackage[T1]{fontenc}
\usepackage{txfonts}
\usepackage{amssymb,verbatim}
\usepackage[all]{xy}
\usepackage{subfig}
\usepackage{tikz}
\usepackage{color}
\usepackage{a4wide}
\usepackage{mathrsfs}
\usepackage{hyperref}
\usepackage{wrapfig}
\usetikzlibrary{arrows}
\DeclareMathAlphabet{\mathpzc}{OT1}{pzc}{m}{it}
\newcommand{\R}{\mathbb{R}}
\newcommand{\C}{\mathbb{C}}

\newcommand\Z{\mathbb{Z}}

\newcommand{\Q}{\mathbb{Q}}

\newcommand{\Pb}{\mathbb{P}}

\newcommand{\Ub}{\mathbf{U}}

\newcommand{\hh}{\mathbf{h}}

\newcommand{\pp}{\mathbf{p}}

\newcommand{\xx}{\mathbf{x}}
\newcommand{\yy}{\mathbf{y}}

\newcommand{\Ccal}{\mathcal{C}}

\newcommand{\FF}{\mathscr{F}}

\newcommand{\OO}{\mathscr{O}}
\newcommand{\PP}{\mathscr{P}}

\newcommand{\Aut}{{\rm Aut}}

\newcommand{\card}{{\rm card}}
\newcommand{\codim}{{\rm codim}}

\newcommand{\id}{\mathrm{id}}
\newcommand{\Id}{\mathrm{Id}}

\newcommand{\CP}{\mathbb{P}_{\mathbb{C}}}
\newcommand{\SL}{{\rm SL}}
\newcommand{\GL}{{\rm GL}}
\newcommand{\Mod}{\mathfrak{M}}

\newcommand{\vol}{{\rm vol}}

\newcommand{\rk}{\mathrm{rk}}

\newcommand{\pr}{\mathrm{pr}}
\renewcommand{\Im}{\mathrm{Im}}
\newcommand{\hg}{\hat{g}}
\newcommand{\hn}{\hat{n}}
\newcommand{\hx}{\hat{x}}
\newcommand{\hX}{\hat{X}}
\newcommand{\hY}{\hat{Y}}
\newcommand{\hy}{\hat{y}}
\newcommand{\hC}{\hat{C}}

\newcommand{\homg}{\hat{\omega}}

\newcommand{\hl}{\hat{\ell}}

\newcommand{\hyy}{\hat{\mathbf{y}}}

\newcommand{\vide}{\varnothing}

\newcommand{\omg}{\omega}
\newcommand{\ol}{\overline}

\newcommand{\cB}{\mathcal{B}}
\newcommand{\cC}{\mathcal{C}}
\newcommand{\cH}{\mathcal{H}}
\newcommand{\cK}{\mathcal{K}}
\newcommand{\cM}{\mathcal{M}}
\newcommand{\cN}{\mathcal{N}}
\newcommand{\cP}{\mathcal{P}}

\newcommand{\cL}{\mathcal{L}}
\newcommand{\cV}{\mathcal{V}}

\newcommand{\cY}{\mathcal{Y}}
\newcommand{\hodge}{\mathcal{H}_{g,n}}
\newcommand{\hodgeext}{\ol{\mathcal{H}}_{g,n}}
\newcommand{\phodge}{\mathbb{P}\mathcal{H}_{g,n}}
\newcommand{\phodgeext}{\mathbb{P}\ol{\mathcal{H}}_{g,n}}
\newcommand{\khodge}{\mathcal{H}^{(k)}_{g,n}}
\newcommand{\khodgeext}{\overline{\mathcal{H}}^{(k)}_{g,n}}
\newcommand{\pkhodge}{\mathbb{P}\mathcal{H}^{(k)}_{g,n}}
\newcommand{\pkhodgeext}{\mathbb{P}\overline{\mathcal{H}}^{(k)}_{g,n}}
\newcommand{\strate}{\cH_{g,n}(\kappa)}

\newcommand{\pstrate}{\mathbb{P}\cH_{g,n}(\kappa)}

\newcommand{\kstrate}{\cH^{(k)}_{g,n}(\kappa)}
\newcommand{\pkstrate}{\mathbb{P}\cH^{(k)}_{g,n}(\kappa)}

\newcommand{\cpkstrate}{\mathbb{P}\ol{\cH}^{(k)}_{g,n}(\kappa)}
\newcommand{\mkstrate}{\cH^{(k)}_{g,\ell}(\kappa)}
\newcommand{\pmkstrate}{\mathbb{P}\cH^{(k)}_{g,\ell}(\kappa)}

\newcommand{\cpmkstrate}{\mathbb{P}\ol{\cH}^{(k)}_{g,\ell}(\kappa)}
\newcommand{\mkhodge}{\cH^{(k)}_{g,\ell}}

\newcommand{\pmkhodgeext}{\mathbb{P}\ol{\cH}^{(k)}_{g,\ell}}
\newcommand{\kstratel}{\cH_{g,\ell}^{(k)}(\kappa)}
\newcommand{\pkstratel}{\mathbb{P}\cH_{g,\ell}^{(k)}(\kappa)}

\newcommand{\cpkstratel}{\mathbb{P}\ol{\cH}^{(k)}_{g,\ell}(\kappa)}
\newcommand{\khodgel}{\cH^{(k)}_{g,\ell}}

\newcommand{\pkhodgel}{\Pb\cH^{(k)}_{g,\ell}}
\newcommand{\pkhodgelext}{\Pb\ol{\cH}^{(k)}_{g,\ell}}
\newcommand{\modac}{\Mod^{\rm ac}_{g,\ell,\kappa}}
\newcommand{\cmodac}{\ol{\Mod}^{\rm ac}_{g,\ell,\kappa}}
\newcommand{\hodgeac}{\mathcal{H}^{\rm ac}_{g,\ell,\kappa}}
\newcommand{\hodgeacext}{\ol{\mathcal{H}}^{\rm ac}_{g,\ell,\kappa}}
\newcommand{\phodgeac}{\mathbb{P}\mathcal{H}^{\rm ac}_{g,\ell,\kappa}}
\newcommand{\phodgeacext}{\mathbb{P}\ol{\mathcal{H}}^{\rm ac}_{g,\ell,\kappa}}
\newcommand{\ihodgeac}{\mathcal{H}^{{\rm ac},\zeta}_{g,\ell,\kappa}}
\newcommand{\ihodgeacext}{\overline{\mathcal{H}}^{{\rm ac},\zeta}_{g,\ell,\kappa}}
\newcommand{\pihodgeac}{\mathbb{P}\mathcal{H}^{{\rm ac},\zeta}_{g,\ell,\kappa}}
\newcommand{\pihodgeacext}{\mathbb{P}\overline{\mathcal{H}}^{{\rm ac},\zeta}_{g,\ell,\kappa}}
\newtheorem{Theorem}{Theorem}[section]
\newtheorem{Corollary}[Theorem]{Corollary}
\newtheorem{Lemma}[Theorem]{Lemma}
\newtheorem{Proposition}[Theorem]{Proposition}

\newtheorem{Claim}[Theorem]{Claim}

\theoremstyle{remark}
\newtheorem{Remark}[Theorem]{Remark}
\begin{document}
\title[On the volumes of linear subvarieties]{On the volumes of linear subvarieties in moduli spaces of projectivized Abelian differentials}
\author{Duc-Manh Nguyen}
\address{Univ. Bordeaux, CNRS, INP Bordeaux, IMB, UMR 5251, F-33405 Talence, France}
\email[D.-M.~Nguyen]{duc-manh.nguyen@math.u-bordeaux.fr}

\date{\today}

\begin{abstract}
For $k \in \Z_{>0}$, let $\mathcal{H}^{(k)}_{g,n}$ denote the vector bundle over $\mathfrak{M}_{g,n}$ whose every fiber consists of  meromorphic $k$-differentials with poles of order at most $k-1$ on a fixed Riemman surface of genus $g$ with $n$ marked points (all the poles must be located at the marked points).
%Denote by $\mathbb{P}\mathcal{H}^{(k)}_{g,n}$ the projective bundle associated to $\mathcal{H}^{(k)}_{g,n}$.
The bundle $\mathcal{H}^{(k)}_{g,n}$ and its associated projective bundle $\mathbb{P}\mathcal{H}^{(k)}_{g,n}$ admit natural extensions, denoted by $\overline{\mathcal{H}}^{(k)}_{g,n}$ and $\mathbb{P}\overline{\mathcal{H}}^{(k)}_{g,n}$ respectively, to the Deligne-Mumford compactification $\overline{\mathfrak{M}}_{g,n}$ of $\mathfrak{M}_{g,n}$.
We prove the following statement: let $\mathcal{M}$ be a subvariety of dimension $d$ of the projective bundle $\mathbb{P}\mathcal{H}^{(k)}_{g,n}$.
Denote by $\mathscr{O}(-1)_{\mathbb{P}\overline{\mathcal{H}}^{(k)}_{g,n}}$ the tautological line bundle over $\mathbb{P}\overline{\mathcal{H}}^{(k)}_{g,n}$.
Then the integral of the $d$-th power of the curvature form of the Hodge norm on $\mathscr{O}(-1)_{\mathbb{P}\overline{\mathcal{H}}^{(k)}_{g,n}}$ over the smooth part of $\mathcal{M}$ is equal to the intersection number of the $d$-th power of the divisor representing $\mathscr{O}(-1)_{\mathbb{P}\overline{\mathcal{H}}^{(k)}_{g,n}}$ and the closure of $\mathcal{M}$ in $\mathbb{P}\overline{\mathcal{H}}^{(k)}_{g,n}$.
As a consequence, if $\mathcal{M}$ is a linear subvariety of the projectivized Hodge bundle $\mathbb{P}\mathcal{H}_{g,n}(=\mathbb{P}\mathcal{H}^{(1)}_{g,n})$ whose local coordinates do not involve  relative periods, then
the volume of $\mathcal{M}$ can be computed  by the self-intersection number of the tautological line bundle on the closure of $\mathcal{M}$ in $\mathbb{P}\overline{\mathcal{H}}_{g,n}(=\mathbb{P}\overline{\mathcal{H}}^{(1)}_{g,n})$.
\end{abstract}
\maketitle
\section{Statements of the main results}\label{sec:intro}
\subsection{Integration of the curvature form of the Hodge norm on moduli spaces of projectivized $k$-differentials}\label{subsec:intro:kdiff:hodge:int}
Let $g,n$ be two non-negative integers such that $2g-2+n>0$.
Let $\ol{\Mod}_{g,n}$ denote the moduli space of complex $n$-pointed stable curves of genus $g$, and $\pi_{g,n}: \ol{\cC}_{g,n}\to \ol{\Mod}_{g,n}$ the universal curve over $\ol{\Mod}_{g,n}$.
As usual, $\Mod_{g,n}$ denotes the subset of $\ol{\Mod}_{g,n}$ which parametrizes  $n$-pointed smooth curves ({\em i.e.} compact Riemann surfaces with $n$ marked points).
Let $K_{\ol{\cC}_{g,n}/\ol{\Mod}_{g,n}}$ denote the relative canonical line bundle of $\pi_{g,n}$, and $D_i$ denote the image of the tautological section $\sigma_i: \ol{\Mod}_{g,n}\to \ol{\cC}_{g,n}$ associated with the $i$-th marked points for $i=1,\dots,n$.

Fix a positive integer $k \in \Z_{>0}$. Define
$$
\cK_{g,n}^{(k)}:=K^{\otimes k}_{\ol{\cC}_{g,n}/\ol{\Mod}_{g,n}}(\sum_{i=1}^n(k-1)\cdot D_i).
$$
For every point $q\in \ol{\Mod}_{g,n}$ denote by $C_q$ the fiber of $\pi_{g,n}$ over $q$. Let $x_i=C_q\cap D_i$ be the $i$-th marked point on $C_q$. Then $\cK^{(k)}_{g,n|C_q} \sim \omg^{\otimes k}_{C_q}(\sum_{i=1}^n(k-1)\cdot x_i)$, where $\omega_{C_q}$ is the dualizing sheaf of $C_q$.
%Note that the restriction of $K_{\ol{\cC}_{g,n}/\ol{\Mod}_{g,n}}$ to $C_q$ is isomorphic to the dualizing sheaf $\omega_{C_q}$ of $C_q$.
By Riemann-Roch Theorem, one readily sees that
$$
\dim H^0(C_q,\cK^{(k)}_{g,n|C_q})= \left\{
\begin{array}{ll}
g & \hbox{ if $k=1$} \\
(g-1)(2k-1)+n(k-1)  & \hbox{ if $k>1$}.
\end{array}
\right.
$$
Since $\dim H^0(C_q,\cK^{(k)}_{g,n|C_q})$ does not depend on $q$, the direct image $\pi_{g,n*}\cK_{g,n}^{(k)}$ is a (holomorphic) vector bundle $\khodgeext$ over $\ol{\Mod}_{g,n}$. The fiber of $\khodgeext$ over $q \in \ol{\Mod}_{g,n}$ is identified with $H^0(C,\omega_{C_q}^{\otimes k}(\sum_{i=1}^n(k-1)\cdot x_i))$, that is the space of $k$-differentials on $C_q$ with poles of order at most $k-1$ at the marked points.
Let $\pkhodgeext$ denote the projective bundle associated with $\khodgeext$.
Elements of $\pkhodgeext$ will be called {\em projectivized $k$-differentials}.

Let $\khodge$ and $\pkhodge$ be the restrictions of $\khodgeext$ and $\pkhodgeext$ to $\Mod_{g,n}$ respectively.
In the case $k=1$, $\cH^{(1)}_{g,n}$ is the usual Hodge bundle. For simplicity, in this case  we will write $\hodgeext, \hodge, \phodgeext, \phodge$ instead of $\ol{\cH}^{(1)}_{g,n}, \cH^{(1)}_{g,n}, \Pb\ol{\cH}^{(1)}_{g,n}, \Pb\cH^{(1)}_{g,n}$ respectively.
Note also that in the case $k=2$, $\cH^{(2)}_{g,n}$ is in fact the cotangent bundle of $\Mod_{g,n}$.

\medskip

By definition, elements of $\khodge\setminus\{0\}$ are tuples $(C,x_1,\dots,x_n,\eta)$, where $C$ is a compact Riemann surface of genus $g$, $x_1,\dots,x_n$ are $n$ marked points on $C$, and $\eta$ is a meromorphic $k$-differential on $C$ such that all of its poles have order at most  $k-1$, and are located in the set $\{x_1,\dots,x_n\}$.  We define
\begin{equation}\label{eq:def:Hodge:norm:kdiff}
||\eta||:=(\int_C|\eta|^{2/k})^{\frac{k}{2}}.
\end{equation}
Since the poles of $\eta$ have orders at most $k-1$, the integral in the right hand side of \eqref{eq:def:Hodge:norm:kdiff} is finite.
Note that $\int_{C}|\eta|^{2/k}$ is the total area of $C$ with respect to the flat metric defined by $|\eta|^{2/k}$.
For this reason, we call $\eta$ a finite area $k$-differential on $C$.
It is not difficult to  check that $||.||$ gives a norm on each fiber of the bundle $\khodge$, which will be called the {\em Hodge norm}.
In the case $k=1$,  $||.||$ is induced by the usual Hodge metric on $\hodge$.
%There are a cyclic ramified cover $f: \hat{C}\to C$ of smooth curves of degree $k$ and an Abelian differential $\homg$ on $\hat{C}$ such that $f^*\eta=\homg^k$.
%The pair $(\hat{C},\homg)$ is called the {\em canonical cyclic cover} of $(C,\eta)$.
%Note that the condition that the poles of $\eta$ have order at most $k-1$ implies that $\homg$ is holomorphic.
%We then define
%$$
%||\eta||:=\left(\frac{\imath}{k}\int_{\hat{C}}\homg\wedge\ol{\homg}\right)^{\frac{k}{2}}
%$$

The Hodge norm provides us with a Hermitian metric on $\OO(-1)_{\phodge}$, the tautological line bundle over $\phodge$.
Let $\Theta$ denote the curvature form of this metric. In this paper we will prove
\begin{Theorem}\label{th:int:curv:form:kdiff}
Let $\cN$ be a subvariety of $\pkhodge$ of dimension $d$. Denote by $\cN_0$ the set of regular points of $\cN$, and by $\ol{\cN}$ the closure of $\cN$ in $\pkhodgeext$. Then we have
\begin{equation}\label{eq:int:cform:self:inters}
\int_{\cN_0}(\frac{\imath}{2\pi}\Theta)^d=c^d_1(\OO(-1)_{\pkhodgeext}) \cdot[\ol{\cN}]
\end{equation}
where $\OO(-1)_{\pkhodgeext}$ is the tautological line bundle over $\pkhodgeext$, and $[\ol{\cN}]$ is the equivalence class of $\ol{\cN}$ in $A_*(\pkhodgeext)$.
\end{Theorem}

If $\ol{\cN}$ is a smooth subvariety of $\pkhodge$ then \eqref{eq:int:cform:self:inters} trivially holds since the cohomology class of $\frac{\imath}{2\pi}\Theta$ is equal to $c_1(\OO(-1)_{\pkhodge})$. The content of Theorem \ref{th:int:curv:form:kdiff} is that \eqref{eq:int:cform:self:inters} also holds when $\ol{\cN}$ is singular and intersects the boundary of $\pkhodgeext$.
In the case $\cN$ is a stratum of $\pkhodge$ (see \textsection \ref{subsec:intro:vol:linear}),  Theorem~\ref{th:int:curv:form:kdiff} has been known by the works \cite{BCGGM:multiscale, CMZ19}.

An unexpected consequence of Theorem~\ref{th:int:curv:form:kdiff} is the following

\begin{Corollary}\label{cor:sub:var:pHodge}
Let $\cN \subset \phodge$ be a subvariety of dimension $d \geq 2g$. Then
$$
c_1^d(\OO(-1)_{\phodgeext})\cdot[\ol{\cN}]=0.
$$
\end{Corollary}
This corollary is in fact not new because it can be derived from the fact that $c_1^{2g}(\OO(-1)_{\phodgeext})= 0 \in H^{4g}(\phodgeext)$, which is a consequence of the Mumford relation $c(E)\cdot c(E^*)=1$. Nevertheless Corrolary~\ref{cor:sub:var:pHodge} gives a stark contrast between the cases $k=1$ and $k>1$. Namely, for $k>1$, let $\cN$ be the principal stratum of $\pkhodge$ (that is the space of $k$-differentials which have poles of order $(k-1)$ at the marked points and only simple zeros elsewhere). In this case $\ol{\cN}=\pkhodgeext$, and $(\imath\Theta)^d$ is a volume form on $\cN$. In particular, the left hand side of \eqref{eq:int:cform:self:inters} must be non-zero. Therefore we have $c^{\dim \pkhodge}_1(\OO(-1)_{\pkhodgeext}) \neq 0 \in H^*(\pkhodgeext)$ if $k>1$.

\subsection{Volumes of linear subvarieties in the moduli spaces of projectivized Abelian differentials}\label{subsec:intro:vol:linear}
A {\em  stratum} of $\khodge$ is a subvariety consisting of (meromorphic) $k$-differentials whose the number of zeros and poles as well as their order are prescribed. Specifically, let $\kappa=(\kappa_1,\dots,\kappa_\ell)$, with $\ell \geq n$, be a sequence of integers such that
\begin{itemize}
\item $\kappa_i \in \Z_{>-k}$, for $i=1,\dots,n$,

\item $\kappa_i \in \Z_{>0}$, for $i=n+1,\dots,\ell$,

\item $\kappa_1+\dots+\kappa_\ell=k(2g-2)$.
\end{itemize}
Then the set
$$
\kstrate:=\{(X,x_1,\dots,x_n,\eta) \in \khodge, \; (\eta)=\kappa_1\cdot x_1+\dots+\kappa_\ell\cdot x_\ell\}
$$
where $x_{n+1},\dots,x_\ell$ are $\ell-n$ distinct points in $X\setminus\{x_1,\dots,x_n\}$ is called a {\em stratum} of $\khodge$.
%It is worth noticing that the all poles of the $k$-differential $\eta$ (if exist) must be contained in the set $\{x_1,\dots,x_n\}$ of the marked points.
By convention, we consider all the points in $\{x_1,\dots,x_n\}$ as zeros or poles of $\eta$, possibly with order $0$.
We will call $\{x_1,\dots,x_\ell\}$ the zero set of $\eta$ and denote it by $Z(\eta)$.

Consider now a stratum $\strate$ of the Hodge bundle $\hodge$. Let $\xx=(X,x_1,\dots,x_n,\eta)$ be a point in $\strate$. Then via {\em period mappings} (see \textsection\ref{subsec:vol:form:lin:subvar} for more details), a neighborhood of $\xx$ in $\strate$ is identified with an open subset of $H^1(X,Z(\eta),\C)$.
Let $\pp: H^1(X,Z(\eta),\C) \to  H^1(X,\C)$ be the natural projection. A {\em linear subvariety} of $\strate$ is an algebraic subvariety $\Omega\cM$ which satisfies the following
\begin{itemize}
\item if $\xx=(X,x_1,\dots,x_n,\eta)\in \Omega\cM$, then the image under the period mappings of any local branch of $\Omega\cM$ through $\xx$ is an open subset of a vector subspace $V$ of $H^1(X,Z(\eta),\C)$, and

\item the restriction of the intersection form on $H^1(X,\C)$ to $\pp(V)$ is non-degenerate.
\end{itemize}
If in addition we have $\ker(\pp)\cap V=\{0\}$, then $\Omega\cM$ is said to be {\em absolutely rigid}.
We denote by $\cM$ the projection of $\Omega\cM$ in $\phodge$.
On an absolutely rigid linear subvariety $\Omega\cM$ of $\strate$, one has a canonical volume form constructed from the intersection form on $H^1(X,\C)$.
By a standard construction (see for instance \cite{Zorich:survey}), this volume form induces in turn a volume form $d\mu$ on $\cM$.

In the case of $k$-differentials with $k\geq 2$, it is a well known  fact that any stratum $\kstrate$ can be locally identified with a linear subvariety of some stratum $\cH_{\hg,\hn}(\hat{\kappa})$ of Abelian differentials via the cyclic cover construction.
For strata $\kstrate$ such that none of the entries of $\kappa$ is divisible by $k$, the corresponding linear subvariety is absolutely rigid.
Thus, for such strata we also have a canonical volume form $d\mu$ on $\pkstrate$ (see \cite{Ng19,KN21}).

\medskip

Among linear subvarieties of $\hodge$, those that are defined locally by linear equations with real coefficients are of particular interest. They are actually orbit closures for an action of $\GL^+(2,\R)$ on $\hodge$, and known as {\em affine invariant submanifolds} of $\hodge$ (see \cite{EM18,EMM15,Filip16}). Their volume is a significant invariant as it allows one to compute other dynamical invariants such as the Lyapunov exponents of the Teichmuller geodesic flow or the Siegel-Veech constants.
It is also well known that the volumes of strata of $k$-differentials, with $k\in \{1,2,3,4,6\}$, give the asymptotics  of the counting of tilings of  surfaces by  triangles and squares (see for instance \cite{EO01, Engel-I, KN21}).

As an application of Theorem~\ref{th:int:curv:form:kdiff} we obtain the following, which is the main motivation of this paper
\begin{Theorem}\label{th:main:vol}\hfill
\begin{itemize}
\item[a)] Let $\Omega\cM$ be an absolutely rigid linear subvariety of $\hodge$ and $\cM$ its projectivization in $\phodge$. Denote by $\ol{\cM}$ the closure of $\cM$ is $\phodgeext$.
Then we have
\begin{equation}\label{eq:vol:lin:subvar}
\mu(\cM)=(-1)^d\cdot\frac{\pi^{d+1}}{(d+1)!}\cdot c_1^{d}(\OO(-1)_{\phodgeext})\cdot [\ol{\cM}]
\end{equation}
where $d=\dim\cM$, and $[\ol{\cM}]$ is the class of $\cM$ in $A_*(\phodgeext)$.

\item[b)] Let $\kstrate$ be a stratum of $\khodge$ and $\pkstrate$ the projectivization of $\kstrate$ in $\pkhodge$ with $k\geq 2$. Let $\cpkstrate$ denote the closure of $\pkstrate$ in $\pkhodgeext$.
Assume that none of the entries of $\kappa$ is divisible by $k$.
Then we have
\begin{equation}\label{eq:vol:str:kdiff}
\mu(\pkstrate)=\frac{(-1)^d}{k^d}\cdot\frac{\pi^{d+1}}{(d+1)!}\cdot c_1^{d}(\OO(-1)_{\pkhodgeext})\cdot [\cpkstrate]
\end{equation}
where $d=\dim\pkstrate$, and $[\cpkstrate]$ is the class of $\cpkstrate$ in $A_*(\pkhodgeext)$.
\end{itemize}
\end{Theorem}
\begin{Remark}\label{rk:vol:coeff}
The constant on the right hand side of \eqref{eq:vol:str:kdiff} depends obviously on the volume form we choose on the stratum $\kstrate$. The values of this constant that correspond to other volume forms on $\kstrate$ have been computed in the literature (see for instance \cite[\textsection 2]{CMS19} and \cite[\textsection 5]{Sau20}).
\end{Remark}

\subsection{Context and related works}\label{subsec:context}
Each stratum of the bundles $\cH^{(1)}_{g,n}$ and $\cH^{(2)}_{g,n}$ carries a special volume form called the {\em Masur-Veech measure}.
The computation of the Masur-Veech volumes of those strata has attracted  great attention because of their application to billiards and Teichm\"uller dynamics.
For a thorough introduction to these fascinating fields of research we refer to the excellent surveys \cite{MT02, Zorich:survey, Wri:survey}.
Masur-Veech measure can also  be defined on strata of $\pkhodge$ for $k\in \{3,4,6\}$ (see for instance \cite{Engel-I}).
In all cases, the Masur-Veech measure always differs from the volume form $d\mu$ in \textsection \ref{subsec:intro:vol:linear} by a constant (see \cite{Ng19} for more details on this constant).

It was shown in \cite{CMSZ20} that the Masur-Veech volumes of the strata of $\phodge$ can be computed by intersection theory on the projectivized Hodge bundle $\phodgeext$.  In \cite{CMZ19}, Costantini, M\"oller, and Zachhuber   showed that \eqref{eq:int:cform:self:inters} holds when $\cM$ is a stratum of $\pkhodge$ for all $k$. Their proof uses in an essential way specific compactifications of strata of $\pkhodge$ known as the {\em moduli spaces of multiscale differentials} that are constructed for  Abelian differentials in \cite{BCGGM:multiscale} and generalized to $k$-differentials for all $k\geq 2$ in \cite{CMZ19}.
The proof we present in this paper does not rely on those constructions.

\medskip

Linear subvarieties of $\hodge$ are generalizations of affine invariant submanifolds  introduced in \cite{EM18} (see also \cite{AEM17} and \cite{Filip16}). The definition of a linear subvariety we use in this paper has been introduced in \cite{KN21}.
In \cite{Beni20} and \cite{BDG20}, Benirschke, Dozier, and Grushevsky  give a similar definition (without the requirement that the restriction of the intersection form to the fiber of the tangent bundle is non-degenerate) and prove some important properties of those subvarieties.
One of the avantages of the definition we use here is that it implies immediately the existence of the volume form $d\mu$. Therefore, the volume with respect to $d\mu$ may be considered as a geometric invariant of absolutely rigid linear subvarieties.

Equality~\eqref{eq:vol:lin:subvar} is an improvement of \cite[Th. 4.1]{KN21} in which we proved a similar equality where on the right hand side instead of $\ol{\cM}$ we have a smooth projective  variety $\hat{\cY}$ admitting a surjective morphism onto $\ol{\cM}$, and where $\OO(-1)_{\phodgeext}$ is replaced by some line bundle on $\hat{\cY}$.

In general, there is no Masur-Veech measure on a linear subvariety unless it is defined by linear equations with rational coefficients,  in which case it is called  {\em arithmetic} (see \cite{Wri:numfield} for more details on the field of definition of an affine invariant submanifold).
The Masur-Veech volume of an arithmetic affine invariant submanifold can be interpreted as the asymptotics of the number of square-tiled surfaces contained in this submanifold. Similarly, the Masur-Veech volumes of strata fo $k$-differentials with $k\in \{2,3,4,6\}$ can also be interpreted as the asymptotics of the numbers of triangulations or quadrangulations on a given (topological) surface.
Theorem~\ref{th:main:vol} means that for absolutely rigid linear subvarities (as well as strata of $\pkhodge$ with $k\in \{2,3,4,6\}$ satisfying the hypothesis in Theorem~\ref{th:main:vol} b)) the Masur-Veech volume can be computed from the self-intersection number of the tautological line bundle modulo the ratio between the Masur-Veech measure and the volume form $d\mu$.

\subsection{Outline of the proof of Theorem~\ref{th:int:curv:form:kdiff}}\label{subsec:outline}
\subsubsection{Case $k=1$.}
The main difficulty of Theorem~\ref{th:int:curv:form:kdiff} in this case is that the Hodge norm does not extend smoothly to the boundary $\partial \phodgeext= \phodgeext\setminus\phodge$ of $\phodgeext$ and that $\ol{\cN}$ may acquire singularities at $\partial \phodgeext$.
Recall that $\Mod_{g,n}$ carries a Variation of Hodge Structure (VHS) of weight 1 $\{F^0,F^1\}$, where $F^0$ is the local system whose fiber over a point $q \in \Mod_{g,n}$ is identified with $H^1(C_q,\C)$, and $F^1 \sim \hodge$.  Note that the local system $F^0$ has unipotent monodromies about the boundary divisor $\partial\ol{\Mod}_{g,n}$ of $\ol{\Mod}_{g,n}$.
%We abusively denote by $\{F^0,F^1\}$ the pullback of this VHS to $\phodge$.

We will show that $\ol{\cN}$ admits a ``resolution'' $\hat{\cN}$ which is in fact an orbifold such that the inverse image of $\partial\ol{\Mod}_{g,n}$ in $\hat{\cN}$ is a normal crossing divisor about which the pullback of the local system $F^0$ has unipotent monodromies. A fundamental result on  VHS  (\cite{Del70, Sch73}) then asserts that $\{F^0,F^1\}$ extends canonically to a filtration of holomorphic vector bundles over $\hat{\cN}$. It follows from  a result of Kawamata (see \cite[p. 266]{Kaw81}) that the canonical extension of $F^1$ is precisely the extended Hodge bundle $\hodgeext$.
Since $\OO(-1)_{\phodgeext}$ is a line subbundle of (the pullback of) $\hodgeext$ over $\phodgeext$, by a deep result by Koll\'ar \cite{Kollar87}, any power of the curvature form of the Hodge norm on $\OO(-1)_{\phodge|\cN}$ is a representative in the sense of currents of the corresponding power of the first Chern class of $\OO(-1)_{\phodgeext}$ on $\hat{\cN}$, and \eqref{eq:int:cform:self:inters} follows.

\subsubsection{Case $k\geq 2$}.
We first show that the problem can be reduced to the case where $\cN$ is subvariety of a stratum $\pkstrate$ of $\pkhodge$, and the (projectivized) $k$-differentials in $\cN$ are primitive, that is they are not tensor powers of some $k'$-differentials with $k' <k$.

One can embed the stratum $\pkstrate$ into $\Pb\cH_{\hg}$ for some  $\hg$ determined by $\kappa$ via the cyclic covering construction. Unfortunately, this embedding does not extend properly to the boundary of $\cpkstrate$. For this reason, we will take a different route  using admissible $\Ub_k$-covers, where $\Ub_k\simeq \Z/k\Z$. To define the appropriate $\Ub_k$-covers, we need to number all the zeros and poles of the $k$-differentials in $\pkstrate$.  This means that we have to pass from $\kstrate \subset \khodge$ to $\mkstrate \subset \mkhodge$, with $\ell=|\kappa|$.  We can then replace $\cN$ by its inverse image in $\pmkstrate$.

Consider now the moduli space $\cmodac$ of admissible $\Ub_k$-covers associated to the stratum $\mkstrate$.  Each element of $\cmodac$ is a stable curve together with a distinguished action of the group $\Ub_k$ such that the quotient is a stable curve in $\ol{\Mod}_{g,\ell}$.
We have a finite morphism $\PP: \cmodac \to \ol{\Mod}_{g,\ell}$ which sends a curve $\hC \in \cmodac$ to $\hC/\Ub_k \in \ol{\Mod}_{g,\ell}$.

Denote by $\hodgeacext$ the Hodge bundle over $\cmodac$ and by $\phodgeacext$  its associated projective bundle.
There is a subbundle $\ihodgeacext$ of $\hodgeacext$, where $\zeta$ is a primitive $k$-th root of unity, together with a finite map $\hat{\PP}:\pihodgeacext \to \pkhodgelext$ which covers the map $\PP: \cmodac \to \ol{\Mod}_{g,\ell}$ and satisfies $\hat{\PP}^*\OO(-1)_{\pkhodgelext} \sim \OO(-1)^{\otimes k}_{\pihodgeacext}$, where as usual $\pihodgeacext$ denote the projective bundle associated to $\ihodgeacext$.
The map $\hat{\PP}$ is  not necessarily  surjective. However, it image contains the closure of $\pmkstrate$. Let $\cM$  be the inverse image of $\cN$ by $\hat{\PP}$. Then equality \eqref{eq:int:cform:self:inters} for $\cN$ is equivalent to the same formula for $\cM$. Now since $\cM$ is a subvariety of the Hodge bundle $\phodgeacext$, the arguments of the case $k=1$ allow us to conclude.

\subsection*{Acknowledgement:}
The author is grateful to Y. Brunebarbe for introducing to him Koll\'ar's works, which actually triggered this project. He warmly thanks V.~ Koziarz and D. Zvonkine for the helpful conversations. The author is partly supported by the French ANR project ANR-19-CE40-0003.

\section{VHS and desingularization}\label{sec:preparation}
\subsection{Extension of the VHS on $\Mod_{g,n}$ to $\ol{\Mod}_{g,n}$}\label{subsec:VHS}
Consider a point $q\in \ol{\Mod}_{g,n}$.
There are an open neighborhood $\widetilde{U}$ of $0$ in $\C^{3g-3}$ (in both Euclidean and Zariski topologies) and  a finite group $G$ acting on $\widetilde{U}$ by isomorphisms  such that  $\widetilde{U}/G$ is isomorphic to a neighborhood of $q$ in $\ol{\Mod}_{g,n}$ with $q$ being identifed with $0$.
Note that in the analytic setting, $G$ can be assumed to be a subgroup of $\GL(3g-3,\C)$.
Moerover, the preimage of $\partial\ol{\Mod}_{g,n}$ in $\widetilde{U}$, which will be denoted by $\partial\widetilde{U}$, is a simple normal crossing divisor in $\widetilde{U}$.
In what follows, we will abusively consider $\widetilde{U}$ as a neighborhood of $q$ in $\ol{\Mod}_{g,n}$.

\medskip

We have a variation of Hodge structure (VHS) of weight 1 $\{F^0,F^1\}$ over $\widetilde{U}^0:=\widetilde{U}\setminus\partial \widetilde{U}$,  where
\begin{itemize}
\item $F^0$ is the complex vector bundle of rank $2g$ associated with the local system whose fiber over a point $(X,x_1,\dots,x_n) \in \widetilde{U}^0$ is given by $H^1(X,\Z)$. By definition, the fiber of $F^0$ over $X$ is identified with
$H^1(X,\Z)\otimes_\Z\C\simeq H^1(X,\C)$.

\item $F^1$ is the holomorphic subbundle of $F^0$ with fiber over $(X,x_1,\dots,x_n)$ being $H^0(X,K_X)$.
\end{itemize}
Note that the total space of $F^1$ is nothing else but $\cH_{g,n|\widetilde{U}^0}$.
The VHS $\{F^0,F^1\}$ comes equipped with the Hodge metric whose  restriction to $F^1$ is precisely $||.||$.

\medskip

Fix a base point  $q_0\sim (X_0,x^0_1,\dots,x^0_n) \in \widetilde{U}^0$. %Let $X_0$ be the Riemann surface represented by $q_0$.
Each irreducible component of $\partial \widetilde{U}$ is determined by the homotopy class of a simple closed curve on the Riemann surface $X_0$ as follows: a generic point in the irreducible component of $\partial\widetilde{U}$ represents a (complex) nodal curve obtained by pinching the corresponding closed curve on $X_0$.
The monodromy of $F^0_{q_0} \simeq H^1(X_0,\C)$ along a loop about an irreducible component of $\partial \widetilde{U}$ is given by a Dehn twist about the corresponding curve on $X_0$. It is a well known fact that the action of a Dehn twist on $H^1(X_0,\Z)$ is given by a transvection matrix in $\mathrm{Sp}(2g,\Z)$ (see for instance \cite[Ch.6]{FM12}). In particular, the associated monodromy is unipotent.
It then follows from a result of Schmid \cite{Sch73} that the filtration $\{F^0,F^1\}$ extends canonically to a filtration $\{\ol{F}^0,\ol{F}^1\}$ of holomorphic vector bundles over $\widetilde U$.

\medskip

Recall that $\pi_{g,n}: \ol{\Ccal}_{g,n} \to \ol{\Mod}_{g,n}$ is the universal curve over $\ol{\Mod}_{g,n}$.
There is a family $\tilde{\pi}_U: \widetilde{\Ccal}_U \to \widetilde{U}$ of $n$-pointed stable curves of  genus $g$ which satisfies
\begin{itemize}
\item the action of $G$ on $\widetilde{U}$ lifts to an action on $\widetilde{\Ccal}_U$ by isomorphisms,

\item for every point $\tilde{q} \in \widetilde{U}$, let $C_{\tilde{q}}$ denote the fiber $\tilde{\pi}_U^{-1}(\{\tilde{q}\})$, then ${\rm Stab}_G(\tilde{p}) \simeq {\rm Aut}(C_{\tilde{p}})$,

\item $\widetilde{\Ccal}_U/G$ is isomorphic to $\pi_{g,n}^{-1}(U) \subset \ol{\Ccal}_{g,n}$, and the restriction of $\pi_{g,n}$ to $\pi^{-1}_{g,n}(U)$ is induced by $\tilde{\pi}_U$.
\end{itemize}
The open $\widetilde{U}$ can be chosen such that $\widetilde{C}_U$  is a smooth variety (that is a complex manifold in analytic setting).
By definition $\tilde{\pi}_U$ has connected fibers, and the restriction of $\tilde{\pi}_U$ to $\tilde{\pi}^{-1}_U(\widetilde{U}^0)$ is a smooth morphism.
Thus it follows from a result by Kawamata (see \cite[Th.5 and \textsection 4]{Kaw81}) that the canonical extension $\ol{F}^1$ of $F^1$ to $\widetilde{U}$ is isomorphic to $\tilde{\pi}_{U*}K_{\widetilde{\Ccal}_U/\widetilde{U}}$, where $K_{\widetilde{\Ccal}_U/\widetilde{U}}$ is the relative canonical bundle of the projection $\tilde{\pi}_U: \widetilde{\Ccal}_U\to \widetilde{U}$.
But by definition, $\tilde{\pi}_{U*}K_{\widetilde{\Ccal}_U/\widetilde{U}}$ is precisely the restriction of  the Hodge bundle $\hodgeext$ to $\widetilde{U}$.
Therefore, we have $\ol{F}^1 \simeq \ol{\cH}_{g,n|\widetilde{U}}$.
Note that the action of $G$ on $\widetilde{U}$ naturally lifts to an action on $\ol{\cH}_{g,n|\widetilde{U}}$ by isomorphisms of vector bundles.

\subsection{Orbifold model}\label{subsec:orbifold}
Consider now a subvariety $\cN$  of $\phodge$. Denote by $\ol{\cN}$ its closure in $\phodgeext$.
Let $h: \phodgeext \to \ol{\Mod}_{g,n}$ be the  natural projection.
Define
$$
\partial \phodgeext:=\phodgeext\setminus \phodge=h^{-1}(\partial \ol{\Mod}_g), \text{ and } \; \partial \ol{\cN}:=\ol{\cN}\setminus\cN.
$$
Note also that $\partial\ol{\cN}=h^{-1}(\partial\ol{\Mod}_{g,n})\cap \ol{\cN}$. Our goal now is to show

\begin{Lemma}\label{lm:orb:desing}
There exist an algebraic orbifold $\hat{\cN}$ and a surjective proper birational morphism $\varphi: \hat{\cN} \to \ol{\cN}$ such that
\begin{itemize}
\item[(i)] $\partial\hat{\cN}:=\varphi^{-1}(\partial\ol{\cN})$ is a normal crossing divisor in $\hat{\cN}$.

\item[(ii)] $\psi:=h\circ \varphi : \hat{\cN}\to \ol{\Mod}_{g,n}$ is an orbifold  morphism.
\end{itemize}
\end{Lemma}
\begin{Remark}\label{rk:orb:desing}\hfill
\begin{itemize}
\item[a)] That $\partial\hat{\cN}$ is a normal crossing divisor in $\hat{\cN}$ means the following: let $p$ be a point in $\partial \hat{\cN}$ with a neighborhood  $V$ isomorphic to a quotient space $\widetilde{V}/H$, where $\widetilde{V}$ is a neighborhood of $0$ in $\C^{\dim \hat{\cN}}$ and $H$ is a finite subgroup of $\GL(\dim \hat{\cN},\C)$. Then the preimage of  $\partial \hat{\cN}$ in $\widetilde{V}$ is defined by an equation of the form $s_1\dots s_a=0$, where $s_1,\dots,s_a$ are some (pairwise distinct) coordinate functions on $\C^{\dim \hat{\cN}}$.

\item[b)] Let $p$ be a point in $\partial\hat{\cN}$ and $q:=\psi(p) \in \ol{\Mod}_{g,n}$. By condition (i) we must have $q \in \partial\ol{\Mod}_{g,n}$.
Let $U=\widetilde{U}/G$ be a neighborhood of $q$ described above.
Condition (ii) means that one can choose $V,U$ such that there exist a group morphism $\rho: H \to G$ and a map $\tilde{\psi}_V:\widetilde{V} \to \widetilde{U}$ satisfying the followings
\begin{itemize}
\item[.] $\tilde{\psi}_V(\gamma\cdot\tilde{p})=\rho(\gamma)\cdot\tilde{\varphi}_V(\tilde{p})$, for all $\tilde{p}\in \widetilde{V}$  and $\gamma\in H$, and

\item[.] the restriction of $\psi$ to $V=\widetilde{V}/H$ is induced by $\tilde{\psi}_V$.
\end{itemize}
Recall that $\partial \widetilde{U}$ is a simple normal crossing divisor in $\widetilde{U}$. Thus there exists a finite family of coordinate functions $t_1,\dots,t_b$ on $\widetilde{U}$ such that $\partial\widetilde{U}$ is defined by the equation $t_1\cdots t_b=0$.
Condition (i) then implies that for all $i\in \{1,\dots,b\}$ we have
$$
\tilde{\psi}^*t_i=u_i\cdot s_1^{\alpha_{i1}}\cdots s_a^{\alpha_{ia}},
\quad \hbox{$\alpha_{ij}\in  \Z_{\geq 0}$},
$$
where $u_i$ is a  non-vanishing holomorphic function on a neighborhood of $p$. This property is crucial to our proof of Therem~\ref{th:int:curv:form:kdiff}.
One can think of $(\hat{\cN},\partial\hat{\cN})$ as a log resolution of $(\ol{\cN},\partial\ol{\cN})$.
However, if we take an arbitrary resolution, then a lift of $\psi_{|V}$ (that is a map from $\widetilde{V}$ to $\widetilde{U}$ with the desired properties) may not exist.
\end{itemize}
\end{Remark}
\begin{proof}[Proof of Lemma~\ref{lm:orb:desing}]
%By definition, $\partial\hat{\cN}$ is a Weil divisor in $\hat{\cN}$.
We first remark that if $\dim\cN=\dim\phodge$, then $\cN=\phodge$. In this case we can simply take $\hat{\cN}=\phodgeext$.
Therefore, from now we will suppose that $\cN$ is a proper subvariety of $\phodge$.

Consider a point $\bar{q}\in \ol{\cN}$, and let $q$ be its projection in $\ol{\Mod}_{g,n}$. Assume first that $q$ is not an orbifold point of $\ol{\Mod}_{g,n}$.
This means that the group $G\simeq \Aut(C_q)$ is trivial, where $C_q$ is the stable curve represented by $q$.
We can choose a neighborhood $U=\widetilde{U}$ of $q$ such that $\ol{\cH}_{g,n|U}$ is trivial.
Hence $\hat{U}:=h^{-1}(U) \subset \phodgeext$ is isomorphic to $U\times \CP^{g-1}$.
Let $\partial\hat{U}$ denote the intersection $\hat{U} \cap \partial\phodgeext$.
Since $\partial U:=\partial \ol{\Mod}_{g,n}\cap U$ is a simple normal crossing divisor in $U$, $\partial \hat{U}$ is also a simple normal crossing divisor in $\hat{U}$.
By definition  $\ol{\cN}_U:=\ol{\cN}\cap \hat{U}$ is a subvariety of $\hat{U}$. Let $I_{\hat{U}}$ denote the ideal sheaf of $\ol{\cN}_U$.

The key ingredient of the proof is an application of the functoriality of the desingularization process (see \cite[Ch. 3]{Kollar07}).
We handle the cases $\codim \cN=1$ and $\codim \cN \geq 2$ in different ways.

\medskip

\noindent \underline{Case $\codim \cN=1$:} in this case $\ol{\cN}_U$ is a divisor in $\hat{U}$. Let $I'_{\hat{U}}$ denote the ideal sheaf of the divisor $\ol{\cN}\cup\partial \hat{U}$ in $\hat{U}$. Apply Theorem 3.35 of \cite{Kollar07} to the triple $(X=\hat{U}, I=I'_{\hat{U}}, E=\varnothing)$, we obtain a sequence of blow-ups
$$
\cB\cP(\hat{U},I'_{\hat{U}},\varnothing)=(\Pi: X_r\overset{\pi_{r-1}}{\to} X_{r-1}\overset{\pi_{r-2}}{\to} \dots \overset{\pi_1}{\to} X_1 \overset{\pi_0}{\to} X_0=\hat{U})
$$
such that the pull-back $\Pi^*I'_U$ is the ideal sheaf of a simple normal crossing divisor $D$ in $X_r$. The proper transform $\ol{\cN}_U^{\#}$  of $\ol{\cN}_U$ in $X_r$ is a finite union $D'_1\cup\dots\cup D'_k$, where each $D'_i$ is an irreducible component of $D$. By definition, each of the $D'_i$ is smooth. However two different components $D'_i, D'_j$ may intersect, therefore $\ol{\cN}^{\#}_U$ is not necessarily smooth.  Let $\hat{\cN}_U$ denote the normalization of $\ol{\cN}_U^{\#}$, which is isomorphic to the disjoint union of the divisor $\{D'_i, \; i=1,\dots,k\}$.
Alternatively, we can define $\hat{\cN}_U$ as the proper transform of $\ol{\cN}^{\#}_U$ after a sequence of blow-ups of $X_r$ along the loci of the intersections of $m$ divisors in the family $\{D'_i, \; i=1,\dots,k\}$, with $m$ taking values in the sequence $(k,k-1,\dots,2)$.
Note that $\hat{\cN}_U$ is smooth, and we have naturally a birational map $\varphi_U: \hat{\cN}_U  \to \ol{\cN}_U$.

We claim that $\varphi^{-1}_U(\partial\hat{U})$ is a simple normal crossing divisor in $\hat{\cN}_U$. By definition, $\Pi^{-1}(\partial \hat{U})$ is  the union of some irreducible components of $D$, none of which is contained in the family $\{D'_i, \; i=1,\dots,k\}$. Since $D$ has simple normal crossing, for each $i \in \{1,\dots,k\}$, the intersection $\Pi^{-1}(\partial \hat{U})\cap D'_i\simeq \varphi^{-1}_{U}(\partial \hat{U})\cap D'_i$ is a simple normal crossing divisor in $D'_i$, and the claim follows.

\medskip

\noindent \underline{Case $\codim\cN \geq 2$:} in this case, we   apply Theorem 3.35 of \cite{Kollar07} to the triple
$$
(X=\hat{U}, I=I_{\hat{U}}, E=\partial\hat{U}).
$$
We then have a sequence of blow-ups
$$
\cB\cP(\hat{U},I_{\hat{U}},\partial\hat{U})=(\Pi: X_r\overset{\pi_{r-1}}{\to} X_{r-1}\overset{\pi_{r-2}}{\to} \dots \overset{\pi_1}{\to} X_1 \overset{\pi_0}{\to} X_0=\hat{U}),
$$
with smooth centers having simple normal crossing with $\partial \hat{U}$ (see \cite[Def. 3.25]{Kollar07}) such that
\begin{itemize}
\item[a)] $\Pi: X_r \to X_0$ is an isomorphism over $X_0\setminus\mathrm{cosupp}(I)=\hat{U}\setminus \ol{\cN}_U$.

\item[b)] The pull-back $\Pi^*I_{\hat{U}}$ is the ideal sheaf of a simple normal crossing divisor.

\item[c)] $\cB\cP$ commutes with smooth morphisms.
\end{itemize}
Since $\codim \cN \geq 2$,  for some $j\in \{1,\dots,r\}$ we have that the center $Z_j\subset X_{j}$ of $\pi_{j}$ contains the proper transform of $\ol{\cN}_U$, but for all $i<j$, the proper transform of $\ol{\cN}_U$ in $X_{i}$ is not contained in the center of the blow-up $\pi_{i}$.
By assumption, $Z_j$ is smooth.
Since $\Pi$ is an isomorphism over $\hat{U}\setminus\ol{\cN}_U$, we must have $\Pi_{j-1}:=\pi_0\circ\dots\circ\pi_{j-1}$ maps $Z_j$ onto  $\ol{\cN}_U$. Hence $\Pi_{j-1|Z_j}: Z_j \to \ol{\cN}_U$ is birational (if $\dim Z_j > \dim \ol{\cN}_U$ then there is a point in $Z_j$ but not in $\Pi_{j-1}^{-1}(\ol{\cN}_U)$, which means that $\ol{\cN}_U \varsubsetneq \Pi_{j-1}(Z_j)$ and there is a point $\bar{q}'$ in $\hat{U}\setminus \ol{\cN}_U$ such that $\Pi^{-1}(\bar{q}')$ has positive dimension).
We claim that $Z_j$ is not contained in $\Pi^{-1}_{j-1}(\partial \hat{U})$.
This is  because $\Pi_{j-1}(Z_j)=\ol{\cN}_U$, while $\partial\hat{U}\cap \ol{\cN}_U$ is a proper subvariety of $\ol{\cN}_U$.
Since $Z_j$ has simple normal crossing with the total transform of $\partial \hat{U}$ in $X_j$ by assumption, we see that $\Pi_{j-1}^{-1}(\partial \hat{U})\cap Z_j$ is a simple normal crossing divisor in $Z_j$.
We can then define $\hat{\cN}_U:=Z_j$.

\medskip

Assume now that $q$ is an orbifold point of $\ol{\Mod}_{g,n}$, that is $G\neq \{1\}$.
Recall that a neighborhood of $q$ in $\ol{\Mod}_{g,n}$ is identified with $U=\widetilde{U}/G$.
Note that the action of $G$ preserves the preimage $\partial\widetilde{U}$ of $\partial\ol{\Mod}_{g,n}$ in $\widetilde{U}$ (this is because the curves corresponding to points in a  $G$-orbit are all isomorphic).
Let $\hat{U}$ denote the pullback  of $\Pb\ol{\cH}_{g,n|U}$ to $\widetilde{U}$.
Then the action of $G$ on $\widetilde{U}$ extends naturally to an action on $\hat{U}$ by bundle isomorphisms.
Note that we can choose $U$ such that $\hat{U}\simeq \widetilde{U}\times\CP^{g-1}$, and  a neighborhood of $\bar{q}$ in $\phodgeext$ can be identified with $\hat{U}/G$.

Let $\widetilde{\cN}_U$ be the preimage of $\ol{\cN}_U$ in $\hat{U}$.
Note that $\widetilde{\cN}_U$ is a $G$-invariant subvariety of $\hat{U}$.
Let $I_{\hat{U}}$ be the ideal sheaf of $\widetilde{\cN}_U$, and $I'_{\hat{U}}$ the ideal sheaf of $\widetilde{\cN}_U\cup\partial \hat{U}$, where $\partial \hat{U}=\hat{U}\cap\partial\phodgeext$.
If $\codim \cN=1$, we apply \cite[Th. 3.35]{Kollar07} to the triple $(\hat{U},I'_{\hat{U}},\varnothing)$. Since $\widetilde{\cN}_U$ and $\partial\hat{U}$ are $G$-invariant, by the functoriality of $\cB\cP(.)$, the proper transform $\widetilde{\cN}^{\#}_U$ of $\widetilde{\cN}_U$ is $G$-invariant. Thus we have a $G$-action on the normalization $\widetilde{\cN}^{\#, {\rm norm}}_U$ of $\widetilde{\cN}^{\#}_U$.
We then define $\hat{\cN}_U:=\widetilde{\cN}^{\#, {\rm norm}}_U/G$.

If $\codim\cN \geq 2$, we apply \cite[Th.3.35]{Kollar07} to the triple $(\hat{U}, I_{\hat{U}},\partial \hat{U})$.
Let $Z_j$ be as above. By the functoriality of the blow-up sequence $\cB\cP(.)$, every element of $G$ lifts to an automorphism of $X_j$ that preserves $Z_j$. In particular, we get an action of $G$ on $Z_j$.
We define $\hat{\cN}_U$ to be the quotient $Z_j/G$.

\medskip

By taking a finite cover of $\ol{\cN}$ by opens of the form $\ol{\cN}\cap (\hat{U}/G)$, with $\hat{U}$ and $G$ as above, we get a finite family of algebraic varieties $\hat{\cN}_U$, each of which is a finite quotient of a smooth one. The functoriality of $\cB\cP(.)$ then allows one to patch these varieties together to form  an algebraic variety $\hat{\cN}$, which is also an orbifold, together with a map $\varphi: \hat{\cN} \to \ol{\cN}$ with the desired properties.
\end{proof}
\begin{Remark}\label{rk:multiscale:desing}
In the case $\cN$ is a stratum  of $\phodge$, the moduli space of (projectivized) multiscale differentials constructed in \cite{BCGGM:multiscale} is an avatar of $\hat{\cN}$.
\end{Remark}

\section{Proof of Theorem~\ref{th:int:curv:form:kdiff} for Abelian differentials}\label{sec:prf:int:curv:form:abel}
We are now in a position to prove the following

\begin{Theorem}\label{th:int:curv:form:abel}
Let $\cN$ be a subvariety of $\phodge$ of dimension $d$. Denote by $\cN_0$ the set of regular points of $\cN$, and by $\ol{\cN}$ the closure of $\cN$ in $\phodgeext$. Let $\Theta$ denote the curvature form of the Hodge norm on $\OO(-1)_{\phodge}$. Then we have
\begin{equation}\label{eq:int:curv:form:abel}
\left(\frac{\imath}{2\pi} \right)^d\int_{\cN_0}\Theta^d=c^d_1(\OO(-1)_{\phodgeext})\cdot[\ol{\cN}]
\end{equation}
where $[\ol{\cN}]$ is the equivalence class of $\ol{\cN}$ in $A_*(\phodgeext)$.
\end{Theorem}

\begin{proof}
Let $\hat{\cN}$ be the complex orbifold obtained by Lemma~\ref{lm:orb:desing}.
By construction, $\hat{\cN}$ has a finite cover by open subsets of the form $\{ V_i=\widetilde{V}_i/G_i, i\in I\}$, where $\widetilde{V}_i$ is a complex manifold and $G_i$ is a finite group acting on $\widetilde{V}_i$ by isomorphisms. Moreover, we have a surjective birational morphism $\varphi: \hat{\cN} \to \ol{\cN}$ such that the composite map $\psi=h\circ \varphi: \hat{\cN} \to \ol{\Mod}_{g,n}$ satisfies the followings: for each $i\in I$, there is an open subset $U_i$ of $\ol{\Mod}_{g,n}$ such that
\begin{itemize}
\item $U_i\simeq \widetilde{U}_i/G_i$, where $\widetilde{U}_i$ is an open neighborhood of $0$ in $\C^{3g-3}$ (in the Euclidean topology), and $G_i$ acts on $\widetilde{U}_i$ by restriction of linear isomorphisms,

\item there is a morphism $\varrho_i: G_i \to \GL(g,\C)$ such that $\ol{\cH}_{g,n|U_i}$ is isomorphic to $\widetilde{U}_i\times \C^g/G_i$, where the action of $G_i$ on $\widetilde{U}_i\times \C^g$ is given by $\gamma\cdot(q,v) = (\gamma\cdot q,\varrho_i(\gamma)\cdot v)$, for all $\gamma \in G_i, q\in \widetilde{U}_i, v \in \C^g$,

%\item the restriction of $\Pb\Omega\ol{\Mod}_{g,n}$ to $U_i$ is isomorphic to $(\widetilde{U}\times\CP^{g-1})/G_i$,

\item $V_i=\psi^{-1}(U_i)$,

\item $\psi_{|V_i}: V_i \to U_i$ lifts to a map $\tilde{\psi}_i: \widetilde{V}_i \to \widetilde{U}_i$  which is $G_i$-equivariant and satisfies $\partial\widetilde{V}_i:=\tilde{\psi}_i^{-1}(\partial \widetilde{U}_i)$ is a  simple normal crossing divisor in $\widetilde{V}_i$, where $\partial\widetilde{U}_i$ is the pre-image of $\partial\ol{\Mod}_{g,n}$ in $\widetilde{U}_i$.
\end{itemize}

Consider a point $p\in \partial\widetilde{V}_i$. By assumption there are some local coordinate functions $s_1,\dots,s_a$ on a neighborhood of $p$ such that $\partial \widetilde{V}_i$ is locally defined by the equation $s_1\cdots s_a=0$. Let $q=\tilde{\psi}_i(p)$. Then there are some local coordinate functions $t_1,\dots,t_b$ of $\widetilde{U}_i$ in a neighborhood of $q$ such that $\partial\widetilde{U}_i$ is defined by the equation $t_1\cdots t_b=0$. The condition $\tilde{\psi}_i^{-1}(\partial\widetilde{U}_i)=\partial\widetilde{V}_i$ means that for all $j\in \{1,\dots,b\}$ we can write
$$
t_j=u_j\cdot s_1^{\alpha{j_1}}\cdots s_a^{\alpha_{j_a}},
$$
where  $\alpha_{j_k} \in \Z_{\geq 0} \; \text{ for all } k=1,\dots,a$, and $u_j$ is a non-vanishing holomorphic function in a neighborhood of $p$.

Let $\widetilde{V}_i^0:=\widetilde{V}_i\setminus\partial\widetilde{V}_i$. By definition, we have $\tilde{\psi}_i(\widetilde{V}_i^0) \subset \widetilde{U}_i^0$.
Recall that we have a VHS $\{F^0,F^1\}$ of weight $1$ over $\widetilde{U}_i^0$. Pulling back by $\tilde{\psi}_i$, we get the VHS $\{\tilde{\psi}_i^*F^0, \tilde{\psi}_i^*F^1\}$ on $\widetilde{V}_i^0$.
Recall that the monodromy of $F^0$ along a loop about any irreducible component of $\partial\widetilde{U}_i$  is unipotent. It follows that the monodromy  of $\tilde{\psi}_i^*F^0$ along a loop  about any irredible component of $\partial\widetilde{V}_i$  is also given by a unipotent matrix. Thus the filtration $\{\tilde{\psi}^*_iF^0,\tilde{\psi}^*_iF^1\}$ extends canonically to $\widetilde{V}_i$, and this extension is isomorphic to $\{\tilde{\psi}^*_i\ol{F}^0,\tilde{\psi}^*_i\ol{F}^1\}$.

In \textsection\ref{subsec:VHS}, we have seen that $\ol{F}^1 \simeq \ol{\cH}_{g,n|\widetilde{U}_i}$. Thus $\tilde{\psi}_i^*\hodgeext$ is the canonical extension of $\tilde{\psi}_i^* F^1$ to $\widetilde{V}_i$.
By construction, the group $G_i$ acts on  $\ol{F}^0, \ol{F}^1$ by isomorphisms of vector bundles. Since $\tilde{\psi}_i$ is $G_i$-equivariant, $G_i$ also acts on $\tilde{\psi}_i^*\ol{F}^1$ by vector bundle isomorphisms.
It follows that the quotient $\tilde{\psi}_i^*\ol{F}^1/G_i$ is an orbifold vector bundle over $V_i$.
Gluing the $\{V_i, \; i\in I\}$ together we obtain an orbifold vector bundle over $\hat{\cN}$  which will be denoted by $\psi^*\ol{F}^1$.
By construction, we have $\psi^*\ol{F}^1$ is isomorphic to $\psi^*\hodgeext$.

\medskip

Consider now the restriction of $\OO(-1)_{\phodgeext}$ to $\ol{\cN}$, which will be denoted by $\cL$.   %Let $\hat{\cL}:=\varphi^*\cL$.
By definition, $\cL$ is a line subbundle of $h^*\ol{\cH}_{g,n|\ol{\cN}}$ (recall that $h:\phodgeext \to \ol{\Mod}_{g,n}$ is the bundle projection). Thus, $\hat{\cL}:=\varphi^*\cL$ is a line subbundle of $\varphi^*\circ h^*\hodgeext=\psi^*\hodgeext$ on $\hat{\cN}$ (in the sense of orbifold vector bundles).

Let $\hat{\cN}^0=\hat{\cN}\setminus\partial\hat{\cN}$.
Since $\psi(\hat{\cN}^0)\subset \Mod_{g,n}$ and the Hodge norm is well defined on $\hodge=\ol{\cH}_{g,n|\Mod_{g,n}}$, we get a well defined Hermitian metric on $||.||$ on $\hat{\cL}_{|\hat{\cN}^0}$. Denote by $\hat{\Theta}$ the curvature form of this metric.
Since $\hat{\cL}$ is a subbundle of $\psi^*\hodgeext \sim \psi^*\ol{F}^1$, by a result of Koll\'ar (see \cite[Theorem 5.1 and Remark 5.19]{Kollar87}) $\left(\frac{\imath}{2\pi}\hat{\Theta}\right)^m$ is a representative in the sense of currents of $c^m_1(\hat{\cL})$ on $\hat{\cN}$ for all $m\in\Z_{\geq 0}$. In particular,  taking $m=d=\dim \ol{\cN}$ we have
\begin{equation}\label{eq:int:on:desing:orb}
\left(\frac{\imath}{2\pi}\right)^{d}\int_{\hat{\cN}^0}\hat{\Theta}^{d}=\left(\frac{\imath}{2\pi}\right)^{d}\int_{\hat{\cN}}\hat{\Theta}^{d}=c^d_1(\hat{\cL})\cdot[\hat{\cN}] \in \Q.
\end{equation}
By definition, we  have $\hat{\cL}=\varphi^*\OO(-1)_{\phodgeext|\ol{\cN}}$  and $\hat{\Theta}=\varphi^*\Theta \hbox{ on $\hat{\cN}^0$}$.
Thus
\begin{equation}\label{eq:int:cform:equal}
\int_{\hat{\cN}^0}\hat{\Theta}^d=\int_{\cN_0}\Theta^d
\end{equation}
and
\begin{equation}\label{eq:self:inters:nb:equal}
c_1^d(\hat{\cL})\cdot[\hat{\cN}]=c_1^d(\OO(-1)_{\phodgeext})\cdot[\ol{\cN}]
\end{equation}
(because $\varphi:\hat{\cN} \to \ol{\cN}$ is birational).
Combining \eqref{eq:int:on:desing:orb}, \eqref{eq:int:cform:equal}, and \eqref{eq:self:inters:nb:equal} we get the desired conclusion.
\end{proof}

%
%Consider a $V_i$ in the family $\{V_i, \; i\in I\}$ above.
%Since objects  on $V_i$ can be viewed as objects defined on $\widetilde{V}_i$ that  are $G_i$-invariant,  we will abusively consider $\widetilde{V}_i$ as a subset of $\hat{\cN}$.
%In particular,  we have $\widetilde{V}_i^0=\widetilde{V}_i\cap \hat{\cN}^0$.
%By construction, we have $\tilde{\psi}_i(\widetilde{V}_i^0)\subset \Mod_{g,n}$.
%Since the Hodge norm is well defined on $\hodge$ and $\hat{\cL} \subset \psi^*\hodgeext$, we get a well defined Hermitian metric on
%$||.||$ on $\hat{\cL}_{|\widetilde{V}_i^0}$. Let $\hat{\Theta}_i$ denote the curvature form of this metric. Note that $\hat{\Theta}_i$ is a closed $(1,1)$-form on $\widetilde{V}_i^0$ which is invariant by $G_i$. Since $\psi^*\hodgeext$ is the VHS canonical extension of $\psi^*F^1\simeq \psi^*\hodge$, and $\hat{\cL}$ is a line subbundle of $\psi^*\hodgeext$, it follows from a result by Kollar (see \cite[Theorem 5.1 and Remark 5.19]{Kollar87}) that for all $m\in \Z_{\geq 0}$, $\left(\frac{\imath}{2\pi}\hat{\Theta}_i\right)^m$ is a representative in the sense of current of $(c^1(\hat{\cL}_{|\widetilde{V}_i}))^m$.
%
%Since the $\hat{\Theta}_i$'s are $G_i$-invariant, they can be patched together to give a well-define $(1,1)$-form $\hat{\Theta}$ on $\hat{\cN}^0$.  By construction, for all $m \in \Z_{\geq 0}$, $(\frac{\imath}{2\pi}\Theta)^m$ is a representative in the sense of current of $(c_1(\hat{\cL}))^m$.

\section{Strata of $k$-differentials and moduli spaces of admissible cyclic coverings}\label{sub:preparation:kdiff}
To prove Theorem~\ref{th:int:curv:form:kdiff} in full generality, we will  reduce the general case to the case $k=1$ and use the arguments of Theorem~\ref{th:int:curv:form:abel} to conclude.
The idea is to use canonical cyclic covers of  $k$-differentials.
More precisely, we will  consider the moduli space of  Abelian differentials that are canonical cyclic covers of $k$-differentials in a fixed connected component of some stratum $\kstrate$.
This space is contained in the total space of the Hodge bundle over the moduli space of stable curves endowed with an admissible $\Ub_k$-action, where $\Ub_k\simeq \Z/k\Z$. Note that the topology of the curves in the latter space and the $\Ub_k$-action strongly depend on the  stratum $\kstrate$. For this reason, we need the subvariety $\cN$ to be contained in  some  stratum of $\khodge$.
This condition is however not at all restrictive, since for any $\cN$ we can always find some stratum that intersects $\cN$ in an open dense subset of $\cN$.
\subsection{Mark and unmark the zeros}\label{subsec:marking:zeros}
Consider a stratum $\kstrate$ of $k$-differentials, where $\kappa=(\kappa_1,\dots,\kappa_\ell)$ with $\ell \geq n$.
Let $(X,x_1,\dots,x_n,\eta)$ be a point in $\kstrate$.
By definition, we have
$$
(\eta)=\sum_{i=1}^\ell \kappa_i\cdot x_i.
$$
Note that the  first $n$ points in the support of $(\eta)$, which can be zeros or poles, are the marked points of $X$, hence they come equipped with a fixed labelling. The remaining points in  the support of $(\eta)$ are all zeros  which are not labelled. To contruct the space of admissible covers of $k$-differentials in $\kstrate$ it is important to label all the zeros and poles of $\eta$, that is we would like to consider $\eta$ as an element of $\mkstrate$. Our first task is to show that the problem in $\kstrate$ can be transposed without change into $\mkstrate$.

Let $Z(\eta)$ denote the support of $(\eta)$.
There is a natural map $\FF: \kstratel \to \kstrate$ which consists of forgetting the labeling of the last $(\ell-n)$ zeros in $Z(\eta)$. It is clear that $\FF$ is a finite morphism, whose degree is equal to the cardinality of the group of permutations of $\{n+1,\dots,\ell\}$ preserving the vector $(\kappa_{n+1},\dots,\kappa_\ell)$.

Let $\pkstrate \subset\pkhodge$ and $\pkstratel\subset \pkhodgel$ be the projectivizations of $\kstrate$ and of $\kstratel$ respectively.
Their closures in  $\pkhodgeext$ and in $\pkhodgelext$ are denoted by $\cpkstrate$ and $\cpkstratel$.  The restrictions of the tautological line bundles on $\pkhodgeext$ and on $\pkhodgelext$ to $\cpkstrate$ and $\cpkstratel$ will be denoted by $\OO(-1)_{\cpkstrate}$ and $\OO(-1)_{\cpkstratel}$.
By a slight abuse of notation,  we will also denote by $\FF$ the induced map from $\pkstratel$ onto $\pkstrate$.

\begin{Proposition}\label{prop:ext:proj:marked:str}
The map $\FF$ admits an extension $\ol{\FF}: \cpkstratel \to \cpkstrate$, and  we have
$$
\ol{\FF}^*\OO(-1)_{\cpkstrate} \sim \OO(-1)_{\cpkstratel}.
$$
\end{Proposition}
\begin{proof}
Consider a point $(Y',y'_1,\dots,y'_\ell,\eta') \in \ol{\cH}^{(k)}_{g,\ell}(\kappa)$. Let $Y'_1,\dots,Y'_r$ be the irreducible components of $Y'$.
It follows from the main result of \cite{BCGGM2} that $\eta'$ corresponds to a collection $(\lambda'_1\cdot\eta'_1,\dots,\lambda'_r\cdot\eta'_r)$,  where $\lambda'_i\in \C$ and $\eta'_i$ is a $k$-differrential on $Y'_i$. For each $i$, the support of $(\eta'_i)$ consists of the marked points that are contained in $Y'_i$, that is $\{y'_1,\dots,y'_\ell\}\cap Y'_i$, and the nodes of $Y'_i$.
Let $(Y,y_1,\dots,y_n) \in \ol{\Mod}_{g,n}$ be the image of $(Y',y'_1,\dots,y'_\ell)$ by the forgetful map.
The curve $Y$ is obtained from $Y'$ by contracting unstable components after removing the marked points $\{y'_{n+1},\dots,y'_\ell\}$. In our situation, such a component is isomorphic to $\CP^1$ and contains  either
\begin{itemize}
\item[(i)] one node, one point in $\{y'_1,\dots,y'_n\}$,  and at least one  point in $\{y'_{n+1},\dots,y'_\ell\}$, or

\item[(ii)] two nodes, no point in $\{y'_1,\dots,y'_n\}$, and at least one point in $\{y'_{n+1},\dots,y'_\ell\}$
\end{itemize}
Assume for instance that $Y'_1$ is a component that is contracted by the forgetful map. Since the order of $\eta'_1$ at a point in $\{y'_1,\dots,y'_n\}$ is at least $1-k$, and at point in $\{y'_{n+1},\dots,y'_\ell\}$ is at least $1$, it follows that in both cases, the order of $\eta'_1$ at a node is smaller than $-k$.
Therefore we must have $\lambda'_1=0$, that is $\eta'$ vanishes identically on $Y'_1$.
Thus by assigning to each component of $Y$ the restriction of $\eta'$ to the corresponding component in $Y'$, we get well defined a $k$-differential $\eta$ which belongs to the fiber of $\khodge$ over $(Y,y_1,\dots,y_n)$. This implies that $\FF: \pkstratel \to \pkstrate$ extends to a map $\ol{\FF}: \cpkstratel \to \pkhodgeext$.
Since $\ol{\FF}$ is a proper, and $\FF(\pkstratel)=\pkstrate$, we get that $\ol{\FF}(\cpkstratel)=\cpkstrate$.
It is also clear from the definition that $\ol{\FF}^*\OO(-1)_{\cpkstrate}$ is isomorphic to $\OO(-1)_{\cpkstratel}$.
\end{proof}

\begin{Remark}\label{rk:ext:proj:marked:str}\hfill
\begin{itemize}
\item   Proposition~\ref{prop:ext:proj:marked:str} does not hold if some of the entries in $\{\kappa_{n+1},\dots,\kappa_\ell\}$ are negative.

\item Assume that $\cN$ is a subvariety of $\pkstrate$. Let $\cN'$ be its pre-image in $\pmkstrate$. Then it follows from Proposition~\ref{prop:ext:proj:marked:str} that equality~\eqref{eq:int:cform:self:inters} holds for $\cN$ if and only if it holds for $\cN'$.
\end{itemize}
\end{Remark}

\subsection{Canonical cyclic covers of $k$-differentials}\label{subsec:cyclic:covers}
Let $(X,x_1,\dots,x_\ell,\eta)$ be a $k$-differential in $\mkstrate$.
We now review some key properties of the canonical cyclic cover of $(X,x_1,\dots,x_\ell,\eta)$.
Recall that $Z(\eta)=\{x_1,\dots,x_\ell\}$.
We first observe that $|\eta|^{2/k}$ defines a flat metric on $X$ with conical singularities in $Z(\eta)$.
The cone angle at $x_i$ is equal to $(1+\frac{{\kappa_i}}{k})\cdot 2\pi$.
The holonomies of the metric $|\eta|^{2/k}$ gives a group morphism
$$
\chi: \pi_1(X\setminus Z(\eta)) \to \Ub_k\simeq \{e^{\frac{2\pi\imath}{k}\cdot j}, \; j=0,\dots,k-1\}.
$$
The image of $\chi$ stays unchanged if we replace $(X,x_1,\dots,x_\ell,\eta)$ by a point nearby in the same stratum.
Thus $\Im(\chi)$ is an invariant of  connected components of $\mkstrate$.
If $\Im(\chi)\neq \Ub_k$, then there exists $k' \in \Z_{>0}$, $k'<k$, such that $\Im(\chi)=\Ub_{k'}$ ($\Ub_{k'}$ is the group of $k'$-th roots of unity).
In this case there is a $k'$-differential $\eta'$ on $X$ such that $\eta={\eta'}^{\otimes k/k'}$.
For this reason, when $\Im(\chi)=\Ub_k$,  we will call $\eta$ a {\em primitive $k$-differential}, and the connected component of $\mkstrate$ to which $(X,x_1,\dots,x_\ell,\eta)$ belongs a {\em primitive component}.
For simplicity, in what follows  we will abusively denote by $\mkstrate$ a connected component of the corresponding stratum, which is supposed to be primitive.

Let us fix a primitive $k$-th root of unity $\zeta$.
It is a well known fact (see \cite{BCGGM2, Ng19} for different proofs) that there are a covering $f: \hX \to X$ of degree $k$ ramified over $Z(\eta)$, an automorphism $\tau  \in \Aut(\hX)$ of order $k$, and an Abelian differential $\homg$ on $\hX$ such that
\begin{itemize}
  \item  $X \simeq \hX/\langle \tau \rangle$,

  \item $\tau^*\homg=\zeta\cdot\homg$, and

  \item $f^*\eta=\homg^k$.
\end{itemize}
Moreover, the triple $(\hX,\homg,\tau)$ is unique up to isomorphism.
It is called the {\em canonical cyclic cover} of $(X,x_1,\dots,x_\ell,\eta)$.
In fact if $Z(\homg)$ be the inverse image of $Z(\eta)$ in $\hX$, then  $f_*\pi_1(\hX\setminus Z(\homg))=\ker(\chi) \subset \pi_1(X\setminus Z(\eta))$.

All of the zeros of $\homg$ are contained in $Z(\homg)$, but it may happen that some of the points in $Z(\homg)$ are not zero of $\homg$. We will refer to these points as  zeros of order $0$ of $\homg$.
By construction, we have a $\Ub_k$-action on $(\hX,Z(\homg))$ which is generated by $\tau$.

\subsection{Admissible $\Ub_k$-covers}\label{subsec:adm:covers}
Let  $\hx_1,\dots,\hx_{\hat{\ell}}$ be the points in $Z(\homg)$.
Let $(Y,y_1,\dots,y_\ell)$ be a pointed smooth curve representing a point in $\Mod_{g,\ell}$.
Call $f_Y: (\hY,\hy_1,\dots,\hy_{\hat{\ell}})  \to (Y,y_1,\dots,y_\ell)$ an admissible $\Ub_k$-cover compatible with $\mkstrate$ if there exist two homeomorphisms $\phi: (X,x_1,\dots,x_\ell) \to (Y,y_1,\dots,y_\ell)$ and $\hat{\phi}:(\hX,\hx_1,\dots,\hx_{\hat{\ell}}) \to (\hY,\hy_1,\dots,\hy_{\hat{\ell}})$ such that $\tau_{\hY}:=\hat{\phi}\circ\tau\circ\hat{\phi}^{-1} \in \Aut(\hY)$ and the following diagram is commutative

\begin{center}
\begin{tikzpicture}[scale=0.4]
\node (A) at (0,4)  {$(\hX,\hx_1,\dots,\hx_{\hat{\ell}})$};

\node (B) at (10,4)  {$(\hY,\hy_1,\dots,\hy_{\hat{\ell}})$};

\node (C) at (0,0) {$(X,x_1,\dots,x_\ell)$};

\node (D) at (10,0) {$(Y,y_1,\dots,y_\ell)$};

\path (A) edge[->, >=angle 60, font=\scriptsize] node[above] {$\hat{\phi}$} (B)
      (A) edge[->, >=angle 60, font=\scriptsize] node[left] {$f$} (C)
      (B) edge[->, >=angle 60, font=\scriptsize] node[right] {$f_Y$} (D)
      (C) edge[->, >=angle 60, font=\scriptsize] node[below] {$\phi$} (D);
\end{tikzpicture}
\end{center}
Two admissible $\Ub_k$-covers $f_{Y}: (\hY,\hy_1,\dots,\hy_{\hat{\ell}})  \to (Y,y_1,\dots,y_{\ell})$  and $f_{Y'}: (\hY',\hy'_1,\dots,\hy'_{\hat{\ell}})  \to (Y',y'_1,\dots,y'_{\ell})$ are isomorphic if there exist two isomorphisms $\hat{\varphi}: \hY \to \hY'$ and $\varphi: Y \to Y'$ such that $f_{Y'}\circ \hat{\varphi}=\varphi\circ f_Y$, and the distinguished automorphisms $\tau_{\hY} \in \Aut(\hY), \tau_{\hY'}\in \Aut(\hY')$ are conjugate by $\hat{\varphi}$.

Note that $\hat{\phi}$ is not unique because precomposing $\hat{\phi}$ with any automorpohism $\tau^i \in \Aut(\hX,\hx_1,\dots,\hx_{\hat{\ell}})$ provides us with another homeomorphism with the same properties. This implies in particular that there is no canonical way to label the points in $\{\hy_1,\dots,\hy_{\hat{\ell}}\}$.

%the diagram
%\begin{center}
%\begin{tikzpicture}[scale=0.4]
%\node (A) at (0,4)  {$(\hY,\hy_1,\dots,\hy_{\hat{\ell}})$};
%
%\node (B) at (10,4)  {$(\hY',\hy'_1,\dots,\hy'_{\hat{\ell}})$};
%
%\node (C) at (0,0) {$(Y,y_1,\dots,y_\ell)$};
%
%\node (D) at (10,0) {$(Y',y'_1,\dots,y'_\ell)$};
%
%\path (A) edge[->, >=angle 60, font=\scriptsize] node[above] {$\hat{\varphi}$} (B)
%      (A) edge[->, >=angle 60, font=\scriptsize] node[left] {$f_Y$} (C)
%      (B) edge[->, >=angle 60, font=\scriptsize] node[right] {$f_{Y'}$} (D)
%      (C) edge[->, >=angle 60, font=\scriptsize] node[below] {$\varphi$} (D);
%\end{tikzpicture}
%\end{center}
%is commutative,

The notion of admissible $\Ub_k$-cover can be generalized to stable curves (see \cite[Ch.16, \textsection 5]{ACG2011}), and there is a (coarse) moduli space of admissible $\Ub_k$-covers compatible with $\mkstrate$, which will be denoted by $\cmodac$.
Elements of $\cmodac$ are tuples $(\hY,\hy_{1},\dots,\hy_{\hat{\ell}},\tau_{\hY})$, where $(\hY,\hy_{1}, \dots,\hy_{\hat{\ell}})$ is a pointed stable curve, and $\tau_{\hY}$ is an automorphism of $\hY$ of order $k$ such that
\begin{itemize}
\item the set $\{\hy_1,\dots,\hy_{\hat{\ell}}\}$ consists of $\ell$ $\langle \tau_{\hY} \rangle$-orbits,

\item let $Y:=\hY/\langle \tau_{\hY} \rangle$ and  $\{y_1,\dots,y_\ell\}$ be the image of $\{\hy_{1},\dots,\hy_{\hat{\ell}}\}$ in $Y$, then $(Y,y_1,\dots,y_\ell)$ is a pointed stable curve in $\ol{\Mod}_{g,\ell}$, and

\item when $\hY$ is smooth, the natural projection $f_Y: (\hY,\hy_{1},\dots,\hy_{\hat{\ell}}) \to (Y,y_1,\dots,y_\ell)$ is an admissible $\Ub_k$-cover compatible with $\mkstrate$.
\end{itemize}
Note that when $\hY$ is a nodal curve, the nodes of $\hY$ are also partitioned into orbits of $\Ub_k$, each orbit maps to a node of $Y$.

\subsection{Projection from $\cmodac$ onto $\ol{\Mod}_{g,\ell}$}\label{subsec:proj:adm:covers:to:Mgl}
There is a natural morphism $\PP: \cmodac \to \ol{\Mod}_{g,\ell}$ which maps a tuple $(\hY,\hy_1,\dots,\hy_{\hat{\ell}},\tau_{\hY}) \in \cmodac$ to the pointed curve $(Y,y_1,\dots,y_\ell)\in \ol{\Mod}_{g,\ell}$.
Let ${\modac}$ denote the set of $(\hY,\hy_1,\dots,\hy_{\hat{\ell}},\tau_{\hY}) \in \cmodac$ where $\hY$ is smooth.
Then $\modac$ is an open dense subset of $\cmodac$. Its complement is denoted by $\partial \cmodac$.

Given $\yy=(Y,y_1,\dots,y_\ell) \in \Mod_{g,\ell}$, each admissible $\Ub_k$-cover of $\yy$ is uniquely detemined by the group morphism $\chi\circ\phi^{-1}_*: \pi_1(Y\setminus\{y_1,\dots,y_\ell\}) \to \Ub_k$.
Since the set of such morphisms is  finite (because $\Ub_k$ is finite), we conclude that the set of admissible $\Ub_k$-covers of $\yy$ is finite. This means that $\PP$ has finite degree, and in particular
$$
\dim \cmodac =\dim \ol{\Mod}_{g,\ell}=3g-3+\ell.
$$

Let $\yy=(Y,y_1,\dots,y_\ell)$ now be point in $\ol{\Mod}_{g,\ell}$.
A neighborhood of $\yy$ in $\ol{\Mod}_{g,\ell}$ can be identified with a quotient space $B/H$, where $B$ is an open  neighborhood of $0$ in $\C^{3g-3+\ell}$ and $H$ a finite group acting on $B$ by isomorphisms.
Let $q_1,\dots,q_\delta$ be the nodes of $Y$. We can choose the local coordinates $(t_1,\dots,t_{3g-3+\ell})$ on $B$ such that $\partial \ol{\Mod}_{g,\ell}\cap B$ is defined by $t_1\cdots t_\delta=0$.

Let $\hyy=(\hY,\hy_{1},\dots,\hy_{\hat{\ell}},\tau_{\hY})$ be a point in $\PP^{-1}(\yy)$. A neighborhood of $\hyy$  in $\cmodac$ can also be identified with a quotient $\tilde{B}/\tilde{H}$, where $\tilde{B}$ is an open subset of $\C^{3g-3+\ell}$, and $\tilde{H}$ is a finite group acting on $\tilde{B}$ by isomorphisms.
Shrinking $\tilde{B}$ if necessary, we can assume that $\PP(\tilde{B}/\tilde{H}) \subset B/H$.
Recall that each node of $Y$ corresponds to a $\Ub_k$-orbit of nodes in $\hY$.  Define $r_j:=\frac{k}{\card f_Y^{-1}(q_j)}$.
There exists a system of local coordinates $(s_1,\dots,s_{3g-3+\ell})$ on $\tilde{B}$ such that $\hyy $ is identified with $0$, and $\tilde{B}\cap\partial\cmodac$ is defined by the equation $s_1\cdots s_\delta=0$.
In this setting, $\PP_{|\tilde{B}/\tilde{H}}: \tilde{B}/\tilde{H} \to B/H$ is induced by a map $\tilde{\PP}_B: \tilde{B} \to B$ which satisfies the following (see \cite[p. 535]{ACG2011})
$$
\tilde{\PP}_B^*t_j = \left\{
\begin{array}{ll}
s_j^{r_j}, & \hbox{ for $j=1,\dots,\delta$}\\
s_j, & \hbox{ for $j=\delta+1,\dots,3g-3+\ell$}.
\end{array}
\right.
$$
This means in particular that the lift of the $r_j$-th power of the Dehn twist associated to the node $q_j$ of $Y$ is the product of the Dehn twists associated to the nodes of $\hY$ in the preimage of $q_j$.

\subsection{Hodge bundle over the space of admissible $\Ub_k$-covers}\label{subsec:hodge:adm:covers}
Denote by $\hodgeacext$ vector bundle over $\cmodac$ whose fiber over $\hyy$ is identified with $H^0(\hY,\omg_{\hY})$.
As usual let $\phodgeacext$ denote the projectivization of $\hodgeacext$. We will call $\hodgeacext$ (resp. $\phodgeacext$) the Hodge bundle (resp. the projectivized Hodge bundle) over $\cmodac$.

Let $\zeta$ be a $k$-th root of unity.
For every $\hyy:=(\hY,\hy_1,\dots,\hy_{\hl},\tau_{\hY}) \in \cmodac$, denote by $H^0(\hY,\omg_{\hY})^\zeta$ the eigenspace associated with $\zeta$ of the action of $\tau_{\hY}$ on $H^0(\hY,\omg_{\hY})$. We have
\begin{Lemma}\label{lm:eig:sp:subbundle}
For all $\zeta \in \Ub_k$, there is a holomorphic subbundle $\ihodgeacext$ of $\hodgeacext$ whose fiber over $\hyy$ is $H^0(\hY,\omg_{\hY})^\zeta$.
\end{Lemma}
\begin{proof}
We essentially only need to check that $\dim H^0(\hY,\omg_{\hY})^\zeta$ does not depend on $\hyy$.
Let $r_{\zeta}(\hyy):=\rk(\tau_{\hY}-\zeta\id)$. Then $r_\zeta$ is a lower continuous function. This means that $r_\zeta(\hyy) \leq r_\zeta(\hyy')$ for all $\hyy'$ in a neighborhood of $\hyy$.
Since $\tau_{\hY}$ has  order $k$, we have
$$
\hg=\dim H^0(\hY,\omg_{\hY})=\sum_{\zeta \in \Ub_k} \dim \ker(\tau_{\hY}-\zeta\id)=\sum_{\zeta \in \Ub_k}\dim\ker(\tau_{\hY'}-\zeta\id)
$$
for all $\hyy'=(\hY',\hy'_1,\dots,\hy'_{\hl},\tau_{\hY'})$ in a neighborhood of $\hyy$. Since $\dim\ker(\tau_{\hY}-\zeta\id) \geq \dim\ker(\tau_{\hY'}-\zeta\id)$ for all $\zeta\in \Ub_k$, actually we must have $\dim\ker(\tau_{\hY}-\zeta\id) = \dim\ker(\tau_{\hY'}-\zeta\id)$
for all $\zeta \in \Ub_k$.
\end{proof}

%There is a natural subbudle $\ihodgeac$ of $\hodgeac$ whose fiber over $\hyy$ is the space $\ker(\tau_{\hY}-\zeta\Id)\subset H^0(\hY,\omg_{\hY})$.
We denote by $\pihodgeacext$ the projectivization of $\ihodgeacext$.
For every $\xi \in H^0(\hY,\omg_{\hY})^\zeta$, we have $\tau_{\hY}^*\xi^{\otimes k}=\xi^{\otimes k}$. This means that $\xi^k$ is the pullback of some $k$-differential $\eta$ on $Y=\hY/\Ub_k$. Taking quotient by the $\C^*$-actions, we get a map $\hat{\PP}: \pihodgeacext \to \pmkhodgeext$
such that the following diagram is commutative
\begin{center}
\begin{tikzpicture}[scale=0.4]
\node (A) at (0,4)  {$\pihodgeacext$};

\node (B) at (10,4)  {$\pmkhodgeext$};

\node (C) at (0,0) {$\cmodac$};

\node (D) at (10,0) {$\ol{\Mod}_{g,\ell}$};

\path (A) edge[->, >=angle 60, font=\scriptsize] node[above] {$\hat{\PP}$} (B)
      (A) edge[->, >=angle 60, font=\scriptsize] node[left] {} (C)
      (B) edge[->, >=angle 60, font=\scriptsize] node[right] {} (D)
      (C) edge[->, >=angle 60, font=\scriptsize] node[below] {$\PP$} (D);
\end{tikzpicture}
\end{center}
where the vertical arrows are bundle projections. Note that even though $\PP$ is surjective, $\hat{\PP}$ is not necessarily surjective because there are $k$-differentials on $Y$ that do not admit $k$-th root on $\hY$.
Nevertheless, by the existence  of the canonical cyclic cover, for some value of $\zeta \in \Ub_k$ (which depends on the choice of the automorphisms $\tau_{\hY}$), the image of $\hat{\PP}$ does contain $\pmkstrate$. From now on we will fix such a value of $\zeta$.
Since $\PP$ is a finite morphism, so is $\hat{\PP}$.

\medskip

Recall that $\modac$ is the set of $\hyy=(\hY,\hy_1,\dots,\hy_{\hat{\ell}},\tau_{\hY})$ where $\hY$ is smooth.
Let $\ihodgeac$ denote the restriction of $\ihodgeacext$ to $\modac$, and by  $\pihodgeac$ its projectivization.
By construction,  the Hodge norms on $\mkhodge$ and on $\hodgeac$ satisfies
the following: if $\xi \in \hodgeac$ is a $k$-th root of $\eta \in \mkhodge$, then we have
$$
||\eta||=||\xi||^k.
$$
The following proposition is straightforward from the construction
\begin{Proposition}\label{prop:proj:ac:to:str}
We have
$$
\hat{\PP}^*\OO(-1)_{\pkhodgelext} \sim \OO(-1)^{\otimes k}_{\pihodgeacext},
$$
where $\OO(-1)_{\pihodgeacext}$ is the tautological line bundle over $\pihodgeacext$.
Denote by $\Theta$ and $\widetilde{\Theta}$ the curvature forms of the Hodge norms on $\OO(-1)_{\pkhodgel}$ and $\OO(-1)_{{\pihodgeac}}$ respectively. Then we have
$$
\hat{\PP}^*\Theta=k\cdot\widetilde{\Theta}
$$
\end{Proposition}

\subsection{Proof of Theorem~\ref{th:int:curv:form:kdiff}}\label{subsec:prf:main:kdiff}
\begin{proof}
In what follows, we will call a connected component of a stratum of $\khodge$ or of $\pkhodge$ simply a stratum.
Given a subset $A$ of $\pkhodgeext$, we denote by $A^0$ its intersection with $\pkhodge$.
Without loss of generality, we can assume that $\cN$ is a closed subvariety  of $\pkhodge$, that is $\ol{\cN}^0=\cN$, and $\cN$ is irreducible.
\begin{Claim}\label{clm:reduce:subv:in:str}
There is a stratum $\pkstrate$ such that $\cN\cap\pkstrate$ is an open dense subset of $\cN$.
\end{Claim}
\begin{proof}
Consider a stratum of highest dimension such that $\cN\cap\pkstrate\neq \vide$. We claim that $\cN \subset\cpkstrate^0$.  Indeed, let  $m=\dim \pkstrate$. Consider  $A=\cpkstrate^0$, and
$$
B=\bigcup_{\substack{\dim\Pb\cH^{(k)}_{g,n}(\kappa') \leq m, \, \kappa'\neq \kappa}}\Pb\ol{\cH}^{(k)}_{g,n}(\kappa')^0
$$
Observe  that $A$ and $B$ are closed subsets of $\pkhodge$, and $A\cup B$ is the union of all strata of dimension at most $m$.
Since $\cN$ does not intersect any stratum of dimension greater than $m$, we have $\cN\subset A\cup B$. Thus $\cN=(\cN\cap A)\cup (\cN\cap B)$.
Since the closure of a stratum cannot intersect another stratum of the same dimension or higher, $B$ does not intersect $\pkstrate$. Therefore we have $\cN\cap B\subsetneq  \cN$ (since $\cN\cap\pkstrate \neq \varnothing$).
If $\cN\cap A\subsetneq \cN$, then we have a contradiction with the hypothesis that $\cN$ is irreducible.
Therefore, we must have $\cN\cap A=\cN$. That is $\cN\subset \cpkstrate^0$.
Since $\pkstrate$ is an open dense subset of $\cpkstrate^0$, $\cN\cap\pkstrate$ is an open dense subset in $\cN$.
\end{proof}

Let $\ell:=|\kappa|$. We now consider the stratum $\mkstrate \subset \khodgel$.  Recall that $\mkstrate$ contains the same $k$-differentals as $\kstrate$, but as elements of $\mkstrate$ all the zeros and poles of those differentials are labelled.
In \textsection~\ref{subsec:marking:zeros}, we have shown that the natural  forgetful map $\FF: \pmkstrate \to \pkstrate$ extends to  a map $\ol{\FF}: \cpmkstrate \to \cpkstrate$ which satisties $\ol{\FF}^*\OO(-1)_{\cpkstrate} \sim  \OO(-1)_{\cpmkstrate}$. Note that the map $\FF$ is  finite-to-one, but over the boundary of $\cpkstrate$  fibers of $\ol{\FF}$ may have positive dimension.

Since $\cN \subset \cpkstrate^0$, we have $\ol{\cN} \subset \cpkstrate$.
Let $\cN', \ol{\cN}'$ be  respectively the preimages of $\cN,\ol{\cN}$ in $\cpmkstrate$.
Note that we have $\cN'=\ol{\cN}'\cap\pkhodgel$.
We abusively denote by $\Theta$ the curvature forms of the Hodge norms on both $\OO(-1)_{\pkstrate}$ and $\OO(-1)_{\pkstratel}$.
\begin{Claim}\label{clm:int:curv:form:marked:str}
Let $\cN'_0$ denote the set of regular points of $\cN'$.
Then equality \eqref{eq:int:cform:self:inters} is equivalent to
\begin{equation}\label{eq:int:equal:self:inters:div:l}
\left(\frac{\imath}{2\pi}\right)^d\int_{\cN'_0}\Theta^d = c^d_1(\OO(-1)_{\cpkstratel})\cdot[\ol{\cN}'].
\end{equation}
\end{Claim}
\begin{proof}
Since $\FF=\ol{\FF}_{|\pkstratel}$ is a finite map, so is the restriction  $\ol{\FF}_{|\cN'\cap\pkstratel}: \cN'\cap\pkstratel \to \cN\cap\pkstrate$. It follows that $\ol{\FF}_{|\ol{\cN}'}: \ol{\cN}'\to \ol{\cN}$ has finite degree.

It is clear from the definition that the forgetful map $\mkstrate\to \kstrate$ preserves  the Hodge norm. Therefore, we have $\FF^*\Theta=\Theta$.
It follows that
$$
\int_{\cN'_0}\Theta^d=\int_{\cN'_0}\FF^*\Theta^d= \delta\cdot\int_{\cN_0}\Theta^d,
$$
where $\delta$ is the degree of the map $\ol{\FF}_{|\ol{\cN}'}$.
On the other hand, we have
\begin{align*}
c^d_1(\OO(-1)_{\cpkstratel})\cdot[\ol{\cN}'] & = \ol{\FF}^*c^d_1(\OO(-1)_{\cpkstrate})\cdot[\ol{\cN}'] \\
& = c^d_1(\OO(-1)_{\cpkstrate})\cdot \ol{\FF}_*[\ol{\cN}']\\
& = \delta\cdot c^d_1(\OO(-1)_{\cpkstrate})\cdot [\ol{\cN}]
\end{align*}
and the claim follows.
\end{proof}
Claim~\ref{clm:int:curv:form:marked:str} means that to prove Theorem~\ref{th:int:curv:form:kdiff}, it suffices to consider the case where $\cN$ is contained in a stratum with all the zeros and poles of the $k$-differentials being labelled.
For this reason, from now on we will assume that $\cN$ is a subvariety of some stratum $\pmkstrate$, where $\ell=|\kappa|$, which is primitive.

\medskip

Consider the moduli space $\cmodac$ of admissible $\Ub_k$-covers compatible with $\mkstrate$. In \textsection\ref{subsec:hodge:adm:covers}, we have seen that there is a  subbundle $\pihodgeacext$ of the Hodge bundle $\hodgeacext$ over $\cmodac$, for some primitive $k$-th root of unity $\zeta$, together with a map $\hat{\PP}:\pihodgeacext \to \pkhodgelext$ which satisfies
\begin{itemize}
\item $\hat{\PP}$ is a finite morphism,

\item $\pmkstrate\subset \hat{\PP}(\pihodgeac)$,

\item $\hat{\PP}^*\OO(-1)_{\pkhodgelext} \sim \OO(-1)^{\otimes k}_{\pihodgeacext}$,

\item $\hat{\PP}^*\Theta= k\cdot\widetilde{\Theta}$, where $\Theta$ and $\widetilde{\Theta}$ are the curvature forms of the Hodge norms on $\OO(-1)_{\pkhodgel}$ and  $\OO(-1)_{\pihodgeac}$ respectively.
\end{itemize}

Since $\hat{\PP}$ is proper, we also have $\cpmkstrate \subset \hat{\PP}(\pihodgeacext)$. Let $\cM$ and $\ol{\cM}$ be the preimages of $\cN$ and $\ol{\cN}$ in $\pihodgeacext$. By the same arguments as Claim~\ref{clm:int:curv:form:marked:str}, we have
\begin{Claim}\label{clm:int:cform:hodge:amd:cov}
Let $\cM_0$ denote the set of regular points in  $\cM$. Then equality~\eqref{eq:int:cform:self:inters} is equivalent to
\begin{equation}\label{eq:int:cform:hodge:amd:cov}
\left(\frac{\imath}{2\pi}\right)^d \int_{\cM_0}\widetilde{\Theta}^d=c_1^d(\OO(-1)_{\phodgeacext})\cdot[\ol{\cM}].
\end{equation}
\end{Claim}
Now, $\cM$ is a subvariety of $\phodgeac$ and $\ol{\cM}$ is its closure in $\phodgeacext$. Thus \eqref{eq:int:cform:hodge:amd:cov} follows from the same arguments as the proof of Theorem~\ref{th:int:curv:form:abel}. Theorem~\ref{th:int:curv:form:kdiff} is then proved.
\end{proof}

\section{Proof of Corollary~\ref{cor:sub:var:pHodge}}\label{prf:cor:sub:var:pHodge}
\begin{proof}
We will actually show that $\Theta^d=0$ everywhere in $\phodge$ for all $d\geq 2g$ and use \eqref{eq:int:cform:self:inters} to conclude.
Let $\xx_0$ be a point in $\phodge$. By definition, $\xx_0$ is a tuple $(X_0,x^0_1,\dots,x^0_n,[\omg_0])$, where $X_0$ is a compact Riemann surface of genus $g$, $x^0_1,\dots,x^0_n$ are $n$ marked points on $X_0$, and $[\omg_0]=\C\cdot\omg_0$ is the complex line generated by a holomorphic $1$-form $\omega_0 \in H^0(X,K_X)$.

Let $\{a_1,\dots,a_g,b_1,\dots,b_g\}$ be a symplectic basis of $H_1(X_0,\Z)$. Without loss of generality, we can assume that $\int_{a_1}\omg_0 \neq 0$.
For all $\xx=(X,x_1,\dots,x_n,[\omg])$ in a neighborhood $U$ of $\xx_0$ in $\phodge$, we can also consider $\{a_1,\dots,a_g,b_1,\dots,b_g\}$ as a symplectic basis of $H_1(X,\Z)$.
Since $\int_{a_1}\omg_0\neq 0$, we also have $\int_{a_1}\omg \neq 0$.
Thus we can normalize $\omg$ by setting $\int_{a_1}\omg=1$ for all $\xx \in U$.
The map $\sigma: \xx \mapsto \omg$ is then a section of $\OO(-1)_{\phodge}$ over $U$.
By definition, we have
$$
\Theta=-\partial\ol{\partial}\log ||\sigma||.
$$
Let us write
$$
z_i:=\int_{a_i}\omg, \; i=2,\dots,g, \quad \text{ and } \quad w_j=\int_{b_j}\omg, \; j=1,\dots,g.
$$
Then $z_i$'s and $w_j$'s are holomorphic functions on $U$.
Now, we have
\begin{align*}
||\sigma(\xx)||^2=||\omg||^2 = \frac{\imath}{2}\cdot\int_X\omg\wedge\ol{\omg} &= \frac{\imath}{2}\cdot \sum_{i=1}^g\left(\int_{a_i}\omg\int_{b_i}\ol{\omg} -\int_{a_i}\ol{\omg}\int_{b_i}\omg \right)\\
                             &= \frac{\imath}{2}\cdot(\bar{w}_1-w_1)+\frac{\imath}{2}\cdot\sum_{i=2}^g(z_i\bar{w}_i-\bar{z}_iw_i).
\end{align*}
Let $E_\xx$ denote the complex vector subspace  generated by $\{dz_i,d\bar{z}_i,\; i=2,\dots,g\}$ and $\{dw_j,d\bar{w}_j, \, j=1,\dots,g\}$ in $T^*_\xx\phodge\otimes \C$. Then $\Theta(\xx) \in \Lambda^{1,1}E_\xx$, and hence $\Theta^d(\xx) \in \Lambda^{d,d}E_\xx$.
But since $\dim_\C E_\xx=2(2g-1)=4g-2$, we have $\Lambda^{d,d}E_\xx=\{0\}$ if $d>2g-1$.
Therefore $\Theta^{2g}$ vanishes identically. It follows from \eqref{eq:int:cform:self:inters} that
$$
c_1^d(\OO(-1)_{\phodgeext})\cdot[\ol{\cN}]=0
$$
whenever $d=\dim \cN \geq 2g$.
\end{proof}

\section{Volumes of absolutely rigid linear subvarieties}\label{sec:vol:abs:rig:subvar}
\subsection{Volume forms on absolutely rigid linear subvariety of $\hodge$}\label{subsec:vol:form:lin:subvar}
Let $\xx:=(X,x_1,\dots,x_n,\omg)$ be a point in some stratum $\strate$ of $\hodge$. For every point $\xx'=(X',x'_1,\dots,x'_n,\omg')$ close enough to $\xx$ in $\strate$, one can identify $H^1(X',Z(\omg'),\C)$ with $H^1(X,Z(\omg),\C)$ via a distinguished homeomorphism between the pairs $(X,Z(\omg))$ and $(X',Z(\omg'))$.
Thus one can associate to $\xx'$  the cohomology class of $\omg'$ in $H^1(X',Z(\omg'),\C)\simeq H^1(X,Z(\omg),\C)$.
This correspondence is called  {\em period mapping}.
A classical result due to H. Masur and W. Veech asserts that period mappings form an atlas of $\strate$ with transition maps given by matrices in $\GL(2g-1+\ell,\Z)$, where  $\ell=|Z(\omg)|$.

Recall that an {\em absolutely rigid linear subvariety} of $\strate$ is an algebraic subvariety $\Omega\cM$ which satisfies the followings:   let $\phi$ be a period mapping defined on a neighborhood $\cV$ of $\xx$:
\begin{itemize}
\item[i)] the image of any irreducible component of $\Omega\cM\cap \cV$ through $\xx$ by $\phi$ is an open subset of a linear subspace
$V$ in $H^1(X,Z(\omg),\C)$,

\item[ii)] the restriction of the intersection form on $H^1(X,\C)$ to $\pp(V)$ is non-degenerate, where $\pp: H^1(X,Z(\omg),\C) \to H^1(X,\C)$ is the natural projection, and

\item[iii)] $V\cap \ker\pp=\{0\}$ (or equivalently, $\pp_{|V}: V \to \pp(V)$ is an isomorphism).
\end{itemize}
On such a subvariety  $\Omega\cM$ we have a natural volume form  defined as follows:  let $\vartheta$ denote the skew-symmetric bilinear form on $H^1(X,\C)$ that is the imaginary part of the intersection form.
%By a slight abuse of notation, we will also denote by $\vartheta$ the pullback of this bilinear form to $H^1(X,Z(\omg),\C)$.
Since $\pp_{|V}: V \to \pp(V)$ is an isomorphism, and the restriction of the intersection form to $\pp(V)$ is non-degenerate by assumption, the restriction of $\vartheta$ to $V$ is a symplectic form.
It follows that $\frac{1}{\dim(V)!}(\vartheta)^{\dim(V)}$ gives a volume form on $V$, which does not depend on the choice of the period mapping.
Therefore, we get a well-defined volume form $d\vol$ on $\Omega\cM$ (see~\cite{Ng19} for more details).

Let $\cM \subset \pstrate$ be the projectivization of $\Omega\cM$.
The volume form $d\vol$ induces a volume form $d\mu$ on $\cM$ by the following well known process:
define
$$
\Omega_1\cM:=\{(X,x_1,\dots,x_n,\omg)\in \Omega\cM, ||\omg|| < 1\}.
$$
Denote by $d\vol_1$ the restriction of $d\vol$ to $\Omega_1\cM$.
The volume form $d\mu$ is then  defined to  be the pushforward of $d\vol_1$ under the projection $\Omega_1\cM \to \cM$.
%The total volume of $\cM$ with respect to $\mu$ is an important invariant of $\cM$ for numerous applications (see for instance \cite{MT02,Zorich:survey}).

\medskip

Consider now a $k$-differential form $\xx:=(X,x_1,\dots,x_n,\eta)$ in some stratum $\kstrate$, with $k\geq 2$.
We assume that $\eta$ is primitive.
Let $(\hX,\hx_1,\dots,\hx_{\hat{\ell}},\homg,\tau)$ be the canonical cyclic cover of $(X,x_1,\dots,x_n,\eta)$.
%By definition $(\hX,\hx_1,\dots,\hx_{\hat{\ell}},\homg)$ belongs to some stratum $\cH_{\hg,\hat{\ell}}(\hat{\kappa})$ of Abelian differentials.
The automorphism $\tau$ acts naturally on the space $H^1(\hX,Z(\homg),\C)$ of  relative cohomology, where $Z(\homg)$ is the inverse image of $\{x_1,\dots,x_n\}$ in $\hX$.
There is a primitive $k$-th root of unity $\zeta$ such that a neighborhood of $\xx$ in $\kstrate$ is identified with an open subset of $V_\zeta:=\ker(\Id -\zeta\tau) \subset H^1(\hX,Z(\homg),\C)$ (see \cite{BCGGM2,Ng19}).
Let $\pp: H^1(\hX,Z(\homg),\C) \to H^1(\hX,\C)$ be the natural projection.
%Denote by $\hh$ the intersection form on $H^1(\hX,\C)$.
We then have (see \cite{Ng19})
\begin{itemize}
\item[a)] $\dim(V_\zeta\cap \ker(\pp))$ equals the number of entries in $\kappa$ that are divisible by $k$, and

\item[b)] the restriction of the intersection form on $H^1(\hX,\C)$ to $\pp(V_\zeta)$ is non-degenerate.
\end{itemize}
It follows in particular that $\kstrate$ can be considered locally as a linear subvariety in some stratum $\cH_{\hg,\hat{\ell}}(\hat{\kappa})$ of Abelian differentials. This linear subvariety is absolutely rigid of none of the entries of $\kappa$ is divisible by $k$.
Thus, in the case where $k$ does not divide any entry of $\kappa$, the spaces $\kstrate$ and $\pkstrate$ come equipped with the volume forms $d\vol$ and $d\mu$ defined above.

For $k\in \{2,3,4,6\}$, we have other natural volume forms on $\kstrate$ (resp. on $\pkstrate$) known as the Masur-Veech volumes (see for instance \cite{Engel-I,KN21}).
The total volume of $\pkstrate$ with respect to the Masur-Veech volume is an important invariant in numerous problems (enumerating tilings of surfaces by triangles or squares, Teichmuller dynamics in moduli spaces, etc...).
It was shown in \cite{Ng19} that the Masur-Veech volume form on $\kstrate$ (resp. on $\pkstrate$) always differs from  $d\vol$ (resp. from $d\mu$) by a rational constant.

\medskip

For the proof of Theorem~\ref{th:main:vol}, we first need to determine the ratio  $\frac{d\mu}{\Theta^d}$ in the case $\cM$ is absolutely rigid. To this purpose, consider a $\C$-vector space $V$ of dimension $m+1$, which is endowed with  a Hermitian form  $\hh$  of signature $(p,q)$, where $p\geq 1$ and $p+q=m+1$. Since $\hh$ is non-degererate, the imaginary part  $\vartheta$ of $\hh$ is a symplectic form on $V$.
Denote by $d\vol$ the volume form $\frac{\vartheta^{m+1}}{(m+1)!}$ on $V$.

Let $V^+$ denote the cone of positive vectors in $V$, that is $V^+=\{v\in V, \; \hh(v,v) >0\}$, and $\Pb V^+:=V^+/\C^* \subset \Pb V$ the projectivization of $V^+$.
Let $V_1=\{v\in V, \, 0< \hh(v,v) <1\} \subset V^+$ and denote by $d\vol_1$ the restriction of $d\vol$ to $V_1$.
We define $\mu$ to be the measure on $\Pb V^+$ which is the pushforward of $d\vol_1$ by the projection $\pr: V_1 \to \Pb V^+$. This means that for every open subset $B\subset \Pb V^+$, $\mu(B)=\vol(C(B)\cap V_1)$, where $C(B)$ is the cone  over $B$ in $V^+$.
The Hermitian form $\hh$ provides us with a  metric on the tautological line bundle $\OO(-1)_{\Pb V^+}$ over $\Pb V^+$. Let $\Theta$ denote the  curvature form of this metric.
\begin{Lemma}\label{lm:ratio:vol:forms}
The measure $\mu$ on $\Pb V^+$ is  defined by a volume form $d\mu$ which satisfies
\begin{equation}\label{eq:ratio:vols}
d\mu=(-1)^m\cdot\frac{2\pi}{2^{m+1}(m+1)!}(\imath\Theta)^{m}.
\end{equation}
\end{Lemma}
\begin{proof}[Sketch of proof]
We can identify $V$ with $\C^{m+1}$ and assume that $\hh$ is given by the matrix $\left(\begin{smallmatrix} I_p & 0 \\ 0 & -I_q \end{smallmatrix} \right)$. Let $(z_0,\dots,z_m)$ be the coordinates on $\C^{m+1}$. We have
$$
\vartheta=\frac{\imath}{2}\left(\sum_{i=0}^{p-1} dz_i\wedge d\bar{z}_i -\sum_{i=p}^{m}dz_i\wedge d\bar{z}_i\right).
$$
and
$$
d\vol=\frac{\vartheta^{m+1}}{(m+1)!}=(-1)^q\cdot\left(\frac{\imath}{2}\right)^{m+1}dz_0\wedge d\bar{z}_0\wedge\dots\wedge d z_{m}\wedge d\bar{z}_m.
$$
Let $\xx$ be a point in $\Pb V^+$, since the  group $U(p,q)$ acts transitively on $\Pb V^+$, we can assume that $\xx=[1:0:\dots:0]$.
We identify a neighborhood  $W$ of $\xx$ in $\Pb V^+$ with a neighborhood of $0$ in $\C^m$ via the map $w:=(w_1,\dots,w_m) \mapsto [1:w_1:\dots:w_m]$.
For all $w=(w_1,\dots,w_m) \in W$, define $\sigma(w):=(1,w_1,\dots,w_m)$,  and
$$
\begin{array}{cccc}
\psi: & [0,2\pi]\times \R_{>0}\times W & \to & V^+\\
   & (\theta,t,w) & \mapsto &  e^{\imath\theta}\cdot t \cdot \sigma(w).
\end{array}
$$
Observe that $\sigma$ is a section of $\OO(-1)_{\Pb V^+}$ over $W$, and
$\psi([0,2\pi]\times \R_{>0}\times W)$ is the cone $C(W)$ over $W$.
Let
$$
f(w):=\hh(\sigma(w),\sigma(w))=1+\sum_{i=1}^{p-1}|w_i|^2-\sum_{i=p}^m|w_i|^2.
$$
Consider a measurable set $B \subset W$.  By definition, we have
\begin{align*}
C(B)\cap V_1 & = \psi(\{(\theta,t,w) \in [0,2\pi]\times \R_{>0}\times B, \; \hh(\psi(\theta,t,w),\psi(\theta,t,w)) <1\})\\
 &= \psi(\{(\theta,t,w) \in [0,2\pi]\times \R_{>0}\times B, \; t^2 < \frac{1}{f(w)}\}).
\end{align*}
A quick computation shows $\psi^*dz_0d\bar{z}_0\dots dz_md\bar{z}_m=-2\imath t^{2m+1} dt d\theta dw_1d\bar{w}_1\dots dw_m d\bar{w}_m$. Therefore
\begin{align*}
\mu(B)= \vol(C(B)\cap V_1) &=(-1)^q\cdot\left(\frac{\imath}{2}\right)^{m}\int_{0}^{2\pi}d\theta \int_B \left(\int_0^{\frac{1}{\sqrt{f(w)}}} t^{2m+1}dt\right) dw_1d\bar{w}_1\dots dw_md\bar{w}_m  \\
&=(-1)^q\cdot\frac{\pi}{m+1}\cdot\frac{1}{f(w)^{m+1}}\cdot\left(\frac{\imath}{2}\right)^{m}\int_B dw_1d\bar{w}_1\dots dw_md\bar{w}_m
\end{align*}
which implies that $\mu$ is induced by the volume form
$$
d\mu=(-1)^q\cdot\frac{\pi}{m+1}\cdot\frac{1}{f(w)^{m+1}}\cdot\left(\frac{\imath}{2}\right)^{m} dw_1d\bar{w}_1\dots dw_md\bar{w}_m.
$$
%Note that since $d\vol$ and $V_1$ are invariant under the action of $U(p,q)$k, we derive that $d\mu$ is always invariant under the action of $U(p,q)$ on $\Pb V^+$.
Now, by definition,
$$
\Theta =-\partial\ol{\partial}\log(f(w))= -\frac{\sum_{i=1}^{p-1}dw_i\wedge d\bar{w}_i-\sum_{i=p}^{m}dw_i\wedge d\bar{w}_i}{f(w)}+\frac{\partial f(w)\wedge \ol{\partial} f(w)}{f^2(w)}
$$
hence
\begin{eqnarray*}
\Theta^m & = & \frac{(-1)^m}{f^{m}(w)}\cdot\left(\sum_{i=1}^{p-1}dw_i\wedge d\bar{w}_i-\sum_{i=p}^{m}dw_i\wedge d\bar{w}_i \right)^m + \\
&  & +m \frac{(-1)^{m-1}}{f^{m+1}(w)}\left(\sum_{i=1}^{p-1}dw_i\wedge d\bar{w}_i-\sum_{i=p}^{m}dw_i\wedge d\bar{w}_i \right)^{m-1}\wedge\partial f(w)\wedge \ol{\partial}f(w)\\
& = & \frac{(-1)^mm!}{f^{m+1}(w)}\left(f(w)-\sum_{i=1}^{p-1}|w_i|^2+\sum_{i=p}^m|w_i|^2\right)\cdot (-1)^q\cdot dw_1\wedge d\bar{w}_1\dots dw_m\wedge d\bar{w}_m\\
& = & \frac{(-1)^{m+q}m!}{f^{m+1}(w)}\cdot dw_1\wedge d\bar{w}_1\dots dw_m\wedge d\bar{w}_m.
\end{eqnarray*}
It follows that
\begin{eqnarray*}
\frac{d\mu}{(\imath \Theta)^m} & = & (-1)^m\cdot \frac{2\pi}{2^{m+1}(m+1)!}
\end{eqnarray*}
and the lemma is proved.
\end{proof}
\begin{Remark}\label{rk:ratio:vol:forms}
The power of $-1$ on the right hand side of \eqref{eq:ratio:vols} is not the same as in the formula in \cite[Lem. 3.2]{KN21}. This is because our definition of the volume form $d\vol$ on $V$ is different from the one in \cite{KN21}. In fact the two volume forms only agree up to sign.
The advantage of the choice of $d\vol$ in this paper (that is $d\vol=\frac{\vartheta^{m+1}}{(m+1)!}$) is that the constant on the right hand side of \eqref{eq:ratio:vols} does not depend on the signature of  $\hh$.
\end{Remark}

\subsection{Proof of Theorem\ref{th:main:vol}}
\begin{proof}
Let $\Omega\cM$ be an absolutely rigid linear subvariety of some stratum $ \strate \subset \hodge$, and $\cM$ its projectivization in $\phodge$.
Let $\Theta$ be the curvature form of Hodge norm on $\OO(-1)_{\Pb\hodge}$.
Recall that $d=\dim\cM$.
Denote by $\cL$  the restriction of $\OO(-1)_{\phodgeext}$ to $\ol{\cM}$.
By Lemma~\ref{lm:ratio:vol:forms}, we have
\begin{equation*}
d\mu=(-1)^d\frac{2\pi}{2^{d+1}(d+1)!}(\imath\Theta)^{d}.
\end{equation*}
Therefore \eqref{eq:vol:lin:subvar} is equivalent to
\begin{equation}\label{eq:int:equal:self:inters:div}
\left(\frac{\imath}{2\pi}\right)^d\int_{\cM}\Theta^{d}=c^d_1(\cL)\cdot[\ol{\cM}]
\end{equation}
which is an immediate consequence of Theorem~\ref{th:int:curv:form:kdiff}.

In the case $\cM=\pkstrate$, where none of the entries of $\kappa$ is divisble by $k$, recall from \textsection \ref{subsec:hodge:adm:covers}, we have a  morphism $\ol{\FF}\circ\hat{\PP}: \pihodgeacext \to \cpkstrate$ such that $(\ol{\FF}\circ\hat{\PP})^*\Theta=k\cdot\widetilde{\Theta}$, where $\Theta$ and $\widetilde{\Theta}$ are respectively the curvature forms of the Hodge norms on $\OO(-1)_{\pkhodgeext}$ and $\OO(-1)_{\phodgeacext}$ (see Proposition~\ref{prop:proj:ac:to:str}). Thus, on $\pkstrate$ we have
\begin{equation*}
d\mu=\frac{(-1)^d}{k^d}\cdot\frac{2\pi}{2^{d+1}(d+1)!}(\imath\Theta)^{d}.
\end{equation*}
and \eqref{eq:vol:str:kdiff} also follows from Theorem~\ref{th:int:curv:form:kdiff}.
\end{proof}

\end{document}